\documentclass[11pt]{article}
\usepackage{amsfonts}

\usepackage{graphics}
\usepackage{indentfirst}
\usepackage{cite}
\usepackage{latexsym}
\usepackage{amsmath}
\usepackage{amssymb}
\usepackage[dvips]{epsfig}
\usepackage{amscd}

\newtheorem{theorem}{Theorem}[section]
\newtheorem{remark}{Remark}[section]

\newtheorem{definition}{Definition}[section]
\newtheorem{lemma}[theorem]{Lemma}

\newtheorem{proposition}[theorem]{Proposition}

\newcommand{\n}{\rho}
\newcommand{\rt}{R_T}
\newcommand{\ti}{\tilde}

\renewcommand{\div}{ {\rm div }  }
\newcommand{\nd}{  \rho^\delta  }
\newcommand{\ud}{  u^\delta  }

\newcommand{\na}{\nabla }
\newcommand{\vp}{\varphi }

\newcommand{\pa}{\partial}

\newcommand{\bi}{\bibitem}

\newcommand{\bt}{\begin{theorem}}
\newcommand{\bl}{\begin{lemma}}
\newcommand{\el}{\end{lemma}}
\newcommand{\et}{\end{theorem}}
\newcommand{\ga}{\gamma}

\newcommand{\curl}{{\rm curl} }
\newcommand{\te}{\theta}

\newcommand{\al}{\alpha}
\newcommand{\de}{\delta}
\newcommand{\ve}{\varepsilon}
\newcommand{\la}{\label}

\newcommand{\ka}{\kappa}

\newcommand{\ol}{\overline}

\newcommand{\bn}{\begin{eqnarray}}
\newcommand{\en}{\end{eqnarray}}
\newcommand{\bnn}{\begin{eqnarray*}}
\newcommand{\enn}{\end{eqnarray*}}

\newcommand{\bnnn}{\begin{eqnarray*}}
\newcommand{\ennn}{\end{eqnarray*}}
\newcommand{\ben}{\begin{enumerate}}
\newcommand{\een}{\end{enumerate}}

\newcommand{\ba}{\begin{aligned}}
\newcommand{\ea}{\end{aligned}}
\newcommand{\be}{\begin{equation}}
\newcommand{\ee}{\end{equation}}

\def\O{\mathbb{T}^2}
\def\p{\partial}
\def\norm[#1]#2{\|#2\|_{#1}}

\def\lap{\triangle}
\def\g{\gamma}
\def\lam{\lambda}

\def\o{\omega}
\makeatletter      
\@addtoreset{equation}{section}
\makeatother       

\title{ Existence and Blowup Behavior of  Global Strong  Solutions  to the Two-Dimensional Baratropic Compressible Navier-Stokes System with   Vacuum and Large Initial Data  \thanks{This research   is partially supported  by
National
Science Foundation of China under  grant 10971215.
 Email: xdhuang@amss.ac.cn (X. Huang), ajingli@gmail.com (J. Li).
 }}

\date{}
\author{Xiangdi H{\small UANG}$^{a}$,  Jing L{\small I}$^{b}$  \\[3mm] {\normalsize $^a$ Academy of Mathematics and System Sciences,} \\
{\normalsize Chinese Academy of Sciences, Beijing 100190, P. R. China} \\[2mm]
{\normalsize $^b$ Institute of Applied Mathematics, AMSS,} \\ {\normalsize \&   Hua Loo-Keng Key Laboratory of Mathematics,}\\
{\normalsize  Chinese Academy of Sciences,    Beijing 100190,
P. R. China}
 }

\begin{document}
\maketitle

\begin{abstract}
  For  periodic   initial data with initial density allowed to vanish, we establish the global existence of strong and weak solutions for the two-dimensional compressible Navier-Stokes equations with    no restrictions on the size of   initial data  provided the shear  viscosity   is a positive constant and the  bulk one is $\lambda = \rho^{\beta}$ with $\beta>4/3$.   These results generalize  and improve  the previous ones due to Vaigant-Kazhikhov([Sib. Math. J. (1995),   36(6),  1283-1316]) which requires $\beta>3$.  Moreover, both the time-independent  upper bound of the density and the large-time   behavior of the strong and weak solutions are also obtained.
\end{abstract}


\section{Introduction and main results}
We study   the two-dimensional barotropic  compressible Navier-Stokes equations which read as follows:
\be\la{n1}
\begin{cases} \rho_t + \div(\rho u) = 0,\\
 (\rho u)_t + \div(\rho u\otimes u) + \nabla P = \mu\lap u + \nabla((\mu + \lam)\div u),
\end{cases}
\ee
where   $\rho=\n(x,t)$ and $u=(u_1(x,t),u_2(x,t))$ represent the density and velocity respectively,  and the pressure $P$ is given by
\be\la{n2}
P(\rho) = a\rho^{\gamma},\quad \ga>1.
\ee
We also have the following hypothesis on the shear viscosity   $\mu$ and the bulk one $\lambda$:
\be\la{n3}
\mu = const,\quad \lam(\n)  = b\rho^{\beta},\,\, b>0,\,\, \beta>0.
\ee
 In the sequel, we set $a=b = 1$ without loosing any generality.

 We consider the Cauchy problem with the given initial data $\n_0$ and $m_0,$ which are
periodic with period $1$ in each space direction $x_i, i = 1, 2,$ i.e., functions defined on
$\O = \mathbb{R}^2/\mathbb{Z}^2.$ We require that
\be \la{n4} \n(x,0)=\n_0(x), \quad \n u(x,0)=m_0(x),\quad x\in \O.\ee

There is a huge literature  concerning the theory of strong and weak solutions for the system of the multidimensional compressible Navier-Stokes equations with constant viscosity coefficients. The local existence and uniqueness of
classical solutions are known in \cite{Na,se1}  in the absence of
vacuum and recently, for strong solutions also, in \cite{cho1,
K2, sal} for the case that the initial density need not be positive
and may vanish in open sets. The global classical solutions were
first obtained by Matsumura-Nishida \cite{M1} for initial data
close to a non-vacuum equilibrium in some Sobolev space $H^s.$   Later, Hoff \cite{Ho4}
studied the problem for discontinuous initial data. For the
existence of solutions for large data,  the major breakthrough is due to
Lions \cite{L1} (see also  Feireisl \cite{F1,Fe}), where he obtained
global existence of weak solutions, defined as solutions with
finite energy, when the exponent $\ga$ is suitably large. The
main restriction on initial data is that the initial energy is
finite, so that the   density is allowed to vanish initially.  Recently, Huang-Li-Xin  \cite{hlx1} established the global existence and uniqueness of classical
solutions to the Cauchy problem for the isentropic compressible Navier-Stokes equations in three-dimensional space with smooth initial data which are of small energy  but possibly large oscillations; in particular, the initial density is allowed to vanish, even has compact support. The compatibility conditions on the initial data of  \cite{hlx1} are further relaxed by \cite{hli1,lzz}.

  However,  there are few results regarding global  strong  solvability for equations of multi-dimensional motions of viscous gas  with    no restrictions on the size of   initial data. One of the first ever ones is due to  Vaigant-Kazhikhov \cite{Ka} who obtained a remarkable result which can be stated that   the two-dimensional system  \eqref{n1}-\eqref{n4} admits a unique global strong solution for large  initial data away from vacuum provided $\beta>3.$   Recently, Perepelitsa \cite{Mik} proved the global existence of a weak solution with uniform lower and upper bounds on the density, as well as the decay of the solution to an equilibrium state in a special case that \be \la{pe1}\beta>3,\quad \ga=\beta,\ee when the initial density is away from vacuum. Very recently, under some additional compatibility  conditions on the initial data,  Jiu-Wang-Xin \cite{jwx} considered classical solutions and  removed  the  condition that the initial density should be away from vacuum in  Vaigant-Kazhikhov \cite{Ka} but still  under the  same condition  $\beta>3$ as that in \cite{Ka}.

Before stating the main results, we explain the notations and
conventions used throughout this paper. We denote
\be\la{1.6}\int fdx=\int_{\O}fdx,\quad \bar{f}=\frac{1}{|\O|}\int f dx .\ee For $1\le r\le \infty ,$  we also denote the standard Lebesgue and
  Sobolev spaces as follows:
$$ L^r=L^r(\O), \quad W^{s,r}= W^{s,r}(\O),  \quad H^s= W^{s,2} .  $$

Then, we give the
definition of weak and strong solutions to \eqref{n1}.

\begin{definition} If  $(\rho,u) $ satisfies \eqref{n1} in the  sense of distribution, then $(\rho,u) $ is called a weak solution to \eqref{n1}.

  If, for  a weak solution, all derivatives involved in \eqref{n1} are regular distributions and   equations  \eqref{n1} hold   almost everywhere   in $\O\times (0,T),$ then the solution is called   strong.
\end{definition}

Thus, the first main result concerning the global existence and large-time behavior of  strong solutions can be stated as follows:
\begin{theorem}\la{t2} Assume that \be\la{bet}\beta>4/3,\quad\ga>1,\ee
 and that the initial data $(\n_0,m_0)$ satisfy that for some $q>2,$
  \be\la{1.9}
  \ba 0\le \rho_0\in   W^{1,q}   ,\quad
 u_0 \in  H^1 ,\quad    m_0=\n_0u_0.
   \ea
  \ee
   Then  the problem  \eqref{n1}-\eqref{n4} has a unique global strong solution $(\n,u)$ satisfying \be\la{1.10}\begin{cases}
  \rho\in C([0,T];W^{1,q} ),\quad  \n_t\in L^\infty(0,T;L^2), \\ u\in L^\infty(0,T; H^1) \cap L^{(q+1)/q}(0,T; W^{2,q}), \\ t^{1/2}u\in L^2(0,T; W^{2,q} )   ,\quad  t^{1/2}u_t\in L^2(0,T;H^1),  \\ \n u\in C([0,T];L^2),   \quad \sqrt{\n} u_t\in L^2(\O\times(0,T)),
   \end{cases}\ee   for any $0<   T<\infty.$
   Moreover, if \be\la{1.11} \beta>3/2, \quad 1< \ga< 3(\beta-1) ,\ee there exists a constant $C$ independent of $T $ such that \be \la{1.12}\sup\limits_{0\le t\le T}\|\n(\cdot,t)\|_{L^\infty }\le C,\ee  and the following large-time  behavior holds: \be \la{1.13}\lim\limits_{t\rightarrow \infty}\left(\|\n-\bar\n_0\|_{L^p}+\|\na u\|_{L^p}\right)=0,\ee  for any $p\in [ 1,\infty).$
\end{theorem}

The second result gives the global existence and large-time behavior of weak solutions.
\begin{theorem}\la{t1} Assume that  \eqref{bet} holds
 and that the initial data $(\n_0,m_0)$ satisfy that
  \be\la{v1.9}
  \ba 0\le \rho_0\in  L^\infty   ,\quad
 u_0 \in  H^1 ,\quad    m_0=\n_0u_0.
   \ea
  \ee
   Then  the problem  \eqref{n1}-\eqref{n4} has at least one  weak solution $(\n,u)$ in $\O\times (0,T)$ for any $T\in (0,\infty).$ Moreover, if $\beta$ and $\ga$ satisfy \eqref{1.11}, there exists a constant $C$ independent of $T $ such that both \eqref{1.12} and \eqref{1.13} hold true.
\end{theorem}

Finally, similar to Li-Xin \cite{lx}, we can obtain from (\ref{1.13}) the
following large-time behavior of the gradient of the density for the strong solution obtained in Theorem \ref{t2} when
vacuum states appear initially.
\begin{theorem} \la{th2} Let $ \beta,\ga$ satisfy \eqref{1.11}. In addition to \eqref{1.9}, assume further
that there exists some point $x_0\in \O $ such that $\rho
_0(x_0)=0.$ Then  the unique global strong solution
$(\rho,u)$ to the Cauchy problem  \eqref{n1}-\eqref{n4}
obtained in Theorem \ref{t2} has to blow up as $t\rightarrow
\infty,$ in the sense that  for any $2<r\le q$ with $q$ as in Theorem \ref{t2},
$$\lim\limits_{t\rightarrow \infty}\|\nabla \rho(\cdot,t)
\|_{L^r}=\infty.$$
\end{theorem}

A few remarks are in order:

\begin{remark}   Theorems \ref{t2} and \ref{t1} generalize  and improve  the earlier results due to Vaigant-Kazhikhov \cite{Ka} where they required that  $\beta>3$ and that the initial density is away from vacuum. \end{remark}

\begin{remark}It should be mentioned here that it seems that  $\beta>1$  is the extremal case for the system \eqref{n1}-\eqref{n3} (see \cite{Ka} or Lemma \ref{kq1}). Therefore, it would be interesting to study the problem \eqref{n1}-\eqref{n4}  when $1< \beta\le 4/3.$ This is left for the future.  \end{remark}

\begin{remark} In   Theorem \ref{t2}, the density is allowed to vanish    initially just under the natural compatibility condition   $m_0=\n_0u_0,$  and no more compatibility ones are required. In fact, our methods can be applied to obtain the local well-posedness of strong solutions to the three-dimensional   system \eqref{n1} just under the natural compatibility condition  $m_0=\n_0u_0.$ This will be reported in a forthcoming paper \cite{hlma}.\end{remark}

\begin{remark} With Theorem \ref{t2} at hand, one can  easily  check that similar to  \cite{hli1, lzz}, if $(\n_0,m_0)$ satisfies for some $q>2,$
  \bnn
  \ba 0\le \rho_0\in   W^{2,q}   ,\quad
 u_0 \in  H^2 ,\quad   m_0=\n_0u_0,
   \ea\enn and the following additional compatibility condition:
\bnn - \mu\lap u_0 - \nabla((\mu + \lam(\n_0))\div u_0)+  \nabla P(\n_0)=\n_0^{1/2}g , \enn   with  some $g\in L^2 ,$ the strong solution obtained in Theorem \ref{t2} becomes
a classical one for positive time. See \cite{hli1,jwx,lzz} for details.  \end{remark}

 \begin{remark}  When the initial density is strictly away from vacuum, Perepelitsa \cite{Mik}  also obtained \eqref{1.12} and  \be \lim\limits_{t\rightarrow \infty}\left(\|\n-\bar\n_0\|_{L^\infty}+\|\na u\|_{L^2}\right)=0,\ee   under the stringent condition  \eqref{pe1}. Note that \eqref{pe1} is a particular case of \eqref{1.11} due to the fact that $3(\beta-1)>\beta$ since $\beta>3/2.$ Thus,     Theorems \ref{t2}  and \ref{t1} improve the results  of Perepelitsa \cite{Mik}.   \end{remark}

We now comment on the analysis of this paper. Note that for smooth initial
data away from vacuum,  the local existence and uniqueness of  strong
solutions to the Cauchy problem  (\ref{n1})-(\ref{n4})  have been
established  in \cite{sal,sol}. Thus, to extend the strong
solutions globally in time and allow the density to vanish initially, one needs global a priori estimates,  which is  independent of the lower  bound of the initial density,  on
smooth solutions to (\ref{n1})-(\ref{n4}) in suitable higher norms. Motivated by our recent studies (\cite{hlx}) on
the blow-up criteria of  strong  solutions to
(\ref{n1}),  it turns out
that the key issue in this paper is to derive  the
  upper bound for the density  which is  independent of the lower  bound of the initial density just under the condition $\beta>4/3.$   To do so, first, similar to \cite{L2,Mik},  we  rewrite $\eqref{n1}_2$ as \eqref{key}
 in terms of  a sum of  commutators of Riesz transforms and the operators of multiplication by $u_i$ (see \eqref{a3.42}).   Then, by energy type estimates   and the compensated  compactness analysis
  \cite[Theorem II.1]{coi3},  we  show that $\log (1+\|\na u\|_{L^2})$
 does not exceed a polynomial  function of $\|\n\|_{L^\infty}$ (see \eqref{n9} and \eqref{nn9}). Next, using the $W^{1,p}$-estimate  of the   commutator  due to  Coifman-Meyer \cite{coi2} (see \eqref{2.7}) and the Brezis-Wainger's inequality   (see \eqref{bmo}),  we obtain an estimate on the $L^\infty$ norm of the commutators in terms of $L^\infty$ norm of the density and $\|\na u\|_{L^2}.$
Both estimates lead to the key a priori
estimate on $\|\n\|_{L^\infty}$  which is independent of the lower  bound of the initial density provided $\beta>4/3.$  See Proposition \ref{aupper} and its proof.

 The
next main step is to bound the gradients of the density  just under the natural compatibility condition   $m_0=\n_0u_0.$    We first obtain the spatial weighted mean estimates on the material derivatives of the velocity which is achieved by modifying the basic estimates on the material derivatives of the velocity  due to Hoff \cite{Ho4}. Then, following  \cite{hlx},   the $L^p$-bound of the gradient  of the density can be obtained by solving a logarithm
Gronwall inequality based on a Beale-Kato-Majda type inequality
(see Lemma \ref{le9}),   the a priori estimates we have just
derived and some careful  initial layer analysis;   and  moreover, such a derivation yields simultaneously
also the bound for $L^1(0,T;L^\infty({\O} ))$-norm of the
gradient of the velocity; see Proposition \ref{le5} and its proof.

The rest of the paper is organized as follows: In Section 2, we collect some
elementary facts and inequalities which will be needed in later analysis. Section 3
is devoted to the derivation
of time-independent and  time-dependent  upper    bounds on the density   which  are independent of the lower of the initial density and needed to extend the local solution to all time. Based on the previous estimates, higher-order ones  are established in Section 4. Then finally, the main results,
Theorems \ref{t2}--\ref{th2}, are proved in Section 5.

\section{Preliminaries}

The following well-known  local existence theory, where the initial
density is strictly away from vacuum, can be found in \cite{sal, sol}.

\begin{lemma}   \la{th0} Assume  that
 $(\n_0,m_0 )$ satisfies \be \la{2.1}
 \n_0\in H^2,  \quad u_0 \in H^2, \quad \inf\limits_{x\in\O}\n_0(x) >0 ,\quad m_0=\n_0u_0.\ee Then there are  a small time
$T >0$ and a  constant $C_0>0$ both depending only on $\|\n_0\|_{H^2},\|u_0\|_{H^2},$ and $\inf\limits_{x\in\O}\n_0(x)$ such that there exists a unique strong solution $(\rho , u )$ to the
  problem   \eqref{n1}-\eqref{n4}  in
$\O\times(0,T )$ satisfying \be\la{2.2}\begin{cases}
  \rho\in C([0,T ];H^2 ),\quad  \rho_t\in C([0,T ];H^1),\\ u\in L^2(0,T ; H^3), \quad u_t\in  L^2(0,T ;H^1),  \\ u_t\in  L^2(0,T ;H^2),  \quad u_{tt}\in  L^2((0,T )\times \O),
   \end{cases}  \ee and \be  \la{2.3}\inf\limits_{(x,t)\in \O\times (0,T )}\n(x,t)\ge C_0>0.\ee
 \end{lemma}

\begin{remark}  It should be mentioned that \cite{sal, sol} dealt  with the case that $\lambda=const.$ However, after some slight modifications, their methods can also be applied to the problem \eqref{n1}-\eqref{n4}.\end{remark}

\begin{remark} In \cite{sal,sol},  instead of $\eqref{2.2}_1,$ it was shown that \bnn \n\in L^\infty(0,T ;H^2),\quad \n_t\in L^\infty(0,T ;H^1). \enn
However, one can use \cite[Lemma 2.3]{L2} to derive  $\eqref{2.2}_1$ by   standard arguments(see \cite{cho1} for details). Moreover, one can also obtain $\eqref{2.2}_3$ by standard arguments due to $\eqref{2.2}_1, \eqref{2.2}_2,$ and $\eqref{2.3}.$  \end{remark}

The   following Poincar\'e-Sobolev and  Brezis-Wainger  inequalities will be used frequently.
\begin{lemma}[\cite{la,en,eng}] \la{leo}
    There exists a
    positive constant $C$ depending only on    $\O $
     such that    every function $u\in H^1(\O )$  satisfies for $2< p<\infty,$
       \be\la{lp}     \|u-\bar u  \|_{L^p }\le C  p^{1/2}\|u-\bar u \|_{L^2 }^{2/p}    \|\na  u\|_{L^2 }^{1-2/p},\quad \|u\|_{L^p }\le C  p^{1/2}\|u\|_{L^2 }^{2/p}    \|u\|_{H^1}^{1-2/p}. \ee
  Moreover, for   $q>2,$     there exists some positive constant $C$ depending only on $q$ and  $\O $ such that  every function $v\in   W^{1,q}(\O )$ satisfies
\be\la{bmo}
\norm[L^\infty]{v}\le C  \|\na v\|_{L^2} \ln^{1/2}(e+\|  \na v\|_{L^q})+C\|v\|_{L^2}+C.
\ee
\end{lemma}

The following Poincar\'e type inequality    can be found in \cite[Lemma 3.2]{Fe}.

\begin{lemma}  \la{jla1} Let  $v\in H^1(\O ),$  and let $\n$ be a non-negative function such that
\bnn 0<M_1\le \int_{\O} \n dx,\quad \int_{\O}  \n^\ga dx\le M_2,\enn
with  $\ga>1.$   Then there is a constant $C$ depending solely on $M_1,M_2$ such that
\be \la{2.5}\|v\|_{L^2(\O )}^2\le C\int_{\O}  \n v^2 dx+C\|\na v\|_{L^2(\O)}^2 .\ee
\end{lemma}

Then, we state the following Beale-Kato-Majda type inequality
which was proved in \cite{B1} when $\div u\equiv 0$ and will be used
later to estimate $\|\nabla u\|_{L^\infty}$ and
$\|\nabla\rho\|_{L^p}$.
\begin{lemma}[\cite{B1,hlx}]   \la{le9}  For $2<q<\infty,$ there is a
constant  $C(q)$ such that  the following estimate holds for all
$\na u\in  W^{1,q} (\O),$ \bnn\la{ww7}\ba \|\na u\|_{L^\infty
}&\le C\left(\|{\rm div}u\|_{L^\infty }+ \|{\rm rot} u\|_{L^\infty }
\right)\log(e+\|\na^2 u\|_{L^q })+C\|\na u\|_{L^2 } +C .
\ea\enn
\end{lemma}

Next, let $\lap^{-1}$ denote  the Laplacian inverse with zero mean on $\O$ and $R_i$ be the usual
 {\sc Riesz} transform on $\O: R_i = (-\lap)^{-1/2}\p_i$. Let $\mathcal{H}^1(\O)$ and $\mathcal{BMO}(\O)$ stand for the usual {\sc Hardy} and {\sc BMO} space:
\bnn
\ba
& \mathcal{H}^1 = \{f\in L^1(\O):\norm[\mathcal{H}^1]{f} = \|f\|_{L^1}+\|R_1f\|_{L^1}+\|R_2f\|_{L^1}
 <\infty ,\quad\bar f = 0\}\\
& \mathcal{BMO} = \{f\in L_{loc}^1(\O):\norm[\mathcal{BMO}]{f}<\infty
\}
\ea
\enn
with
\bnn
\norm[\mathcal{BMO}]{f} = \sup_{x\in\O,r\in (0,d)}\frac{1}{|\Omega_r(x)|}
\int_{\Omega_r(x)}\left|f(y)-\frac{1}{|\Omega_r(x)|}\int_{\Omega_r(x)}f(z)dz\right|dy,
\enn
where $d$ is the diameter of $\O,$ $\Omega_r(x) = \O\cap B_r(x),$ and $B_r(x)$ is a ball with
center $x$ and radius $ r. $
Consider the composition of two Riesz transforms, $ R_i \circ R_j  (i, j = 1, 2).$ There is a
representation of this operator as a singular integral
\bnn R_i \circ R_j (f)(x)=p.v. \int K_{ij}(x-y)f(y)dy,\enn
where the kernel $K_{ij}(x) (i, j = 1, 2)$  has a singularity of the second order at $0$ and
$$|K_{ij}(x)|\le C|x|^{-2}, \,\, x\in \O.$$
Given a function $b,$ define the linear operator
$$[b,R_iR_j](f)\triangleq bR_i\circ R_j(f)-R_i\circ R_j (bf), \,\,i,j= 1, 2.$$
This operator can be written as a convolution with the singular kernel $K_{ij},$
$$[b,R_iR_j](f)(x)\triangleq p.v. \int K_{ij}(x-y)(b(x)-b(y))f(y)dy,\,\,i,j= 1, 2. $$

The following properties of the  commutator $[b,R_iR_j](f)$   will be   useful for our analysis. \begin{lemma} Let $b,f\in C^\infty(\O).$  Then for $p\in (1,\infty),$ there is $ C(p)$ such that
\be\la{2.6} \|[b,R_iR_j](f)\|_{L^p}\le C(p)\|b\|_{\mathcal{BMO}}\|f\|_{L^p}.
\ee
Moreover, for $q_i\in (1,\infty) (i=1,2,3)$ with $q_1^{-1}=q_2^{-1} +q_3^{-1},$ there is $ C$ depending only on $q_i (i=1,2,3)$ such that
\be \la{2.7}\|\na[b,R_iR_j](f)\|_{L^{q_1}}\le C \|\na b\|_{L^{q_2}}\|f\|_{L^{q_3}}.
\ee
\end{lemma}

\begin{remark} Properties \eqref{2.6} and \eqref{2.7} are  due to Coifman-Rochberg-Weiss \cite{coi1} and Coifman-Meyer \cite{coi2} respectively.
\end{remark}

Finally, the following Zlotnik  inequality will be used to get the
uniform (in time) upper bound of the density $\n.$
\begin{lemma}[\cite{zl1}]\la{ale1}   Let the function $y$ satisfy
\bnn y'(t)= g(y)+h'(t) \mbox{  on  } [0,T] ,\quad y(0)=y^0, \enn
with $ g\in C(R)$ and $y,h\in W^{1,1}(0,T).$ If $g(\infty)=-\infty$
and \be\la{a100} h(t_2) -h(t_1) \le N_0 +N_1(t_2-t_1)\ee for all
$0\le t_1<t_2\le T$
  with some $N_0\ge 0$ and $N_1\ge 0,$ then
\bnn y(t)\le \max\left\{y^0,\ti{\zeta} \right\}+N_0<\infty
\mbox{ on
 } [0,T],
\enn
  where $\ti{\zeta} $ is a constant such
that \bnn\la{a101} g(\zeta)\le -N_1 \quad\mbox{ for }\quad \zeta\ge \ti{\zeta}.\enn
\end{lemma}

\section{\la{se3}A priori estimates (I): upper bound of the density}

In this section and the next, we will always assume that $(\n_0,m_0)$ satisfies \eqref{2.1} and  $(\rho,u)$ is the strong solution to  \eqref{n1}-\eqref{n4} on $\O\times (0,T]$ obtained by Lemma \ref{th0}.

\subsection{Time-independent upper bound of the density}

In this subsection, we will establish the following time-independent upper bound of the  density provided \eqref{1.11} holds. Throughout this subsection, we   use the convention that $C$ denotes a
generic positive constant
  independent of both the time $T$ and the lower bound of the initial  density,  and   we write $C(\al)$ to emphasize
that $C$ depends on $\al.$
  \begin{proposition}\la{lle5}   If   \eqref{1.11} holds,   there is a positive constant  $C$ depending only on $  \mu,  \beta, \gamma,  \|\n_0\|_{L^\infty},$ and $\| u_0\|_{H^1} $  such that
  \be\la{a3.56}
   \sup\limits_{0\le t\le T}\|\n\|_{L^\infty}\le C.
  \ee
\end{proposition}

Before proving Proposition \ref{lle5}, we establish a series of a priori estimates, Lemmas \ref{le1}-\ref{lce}.
To proceed, we denote by \bnn \nabla^{\perp} = (\p_2,-\p_1),\quad \frac{D}{Dt}f  =\dot{f} = f_t  + u \cdot\nabla f ,\enn where $\frac{D}{Dt}f$ is the material derivative of $f.$ Let $G$ and $\o$ be the effective viscous flux and the vorticity respectively as follows:
\bnn\la{nn4}
G \triangleq (2\mu+\lam(\rho))\div u - (P-\bar P),\quad \o \triangleq \nabla^{\perp}\cdot u= \p_2u_1 - \p_1u_2.
\enn
Then, we define
\be\la{an7}
( A_1(t))^2 \triangleq \int_{\O}\left((\o(t))^2 + \frac{(G(t))^2}{2\mu+\lam(\n(t))}\right)dx,\ee \be
  \la{n7}(A_2(t))^2  \triangleq \int_{\O}\rho(t)|\dot{u}(t)|^2dx,\ee \be
 \la{bn7} (A_3(t))^2  \triangleq \int_{\O}\left((2\mu+\lam(\n(t)))(\div u(t))^2 + \mu(\o(t))^2\right)dx,
 \ee
 and \be  \la{cn7} \rt \triangleq\sup\limits_{0\le t\le T} \|\n\|_{L^\infty}.\ee
Without loss of generality, we assume that \bnn \la{wq1}
\int \n_0dx=1,
\enn
which together with   $(\ref{n1})_1$  gives
\be  \la{c3.10}  \rt\ge  \|\n(\cdot,t)\|_{L^\infty}\ge \int \n(x,t) dx =\int \n_0dx =1.\ee

Then, we have the following lemma.
\begin{lemma}\la{le1} For any $\al\in(0,1),$   there is a positive constant  $C(\al)$ depending  only on $ \al, \mu,  \beta, \gamma,  \|\n_0\|_{L^\infty},$ and $\| u_0\|_{H^1} $ such that
  \be\la{a3.14}
  \ba
   \frac{d}{dt}A_1^2(t) + A_2^2(t) & \le C(\al) \left(\rt  \vp^2 +  \|\n\|_{L^\beta}^{\beta/2}\vp \right)A_3^2,
  \ea
  \ee   where $\vp$ is defined by
\be\la{a3.15}\ba
\vp(t)\triangleq   1+ A_1 \rt^{ \al\beta /2 }    +   \left\|\frac{ P }{(2\mu+\lam)^{3/2}}\right\|_{L^4}+ \left\| \frac{P }{(2\mu+\lambda)^{1/2}}\right\|_{L^2} . \ea
\ee

\end{lemma}

{ \it Proof. }  First,
the standard energy inequality reads:
\be\la{r3}
 \sup\limits_{0\le t\le T}\int\left(\n|u|^2+\n^\ga\right)dx +\int_0^T A_3^2(t)dt \le C,
\ee which together with \eqref{2.5}  gives that  for $t\in [0,T],$
\be \la{a3.11}
  CA^2_3(t)- C-C\rt^{\ga-\beta} \le A^2_1(t) \le CA^2_3(t) +  C+C\rt^{\ga-\beta} ,  \ee
   \be \la{a3.12}
  C^{-1}\|\na u\|_{L^2}^2 \le A_3^2(t)\le C\rt^{ {\beta} }\|\na u\|_{L^2}^2,  \ee
and \be\la{u-1} \|u\|_{H^1}\le C+C\|\na u\|_{L^2}\le C+CA_3.\ee

Next,   direct calculations show that
\be\la{ra3}\ba
  \nabla^{\perp}\cdot \dot u = \frac{D}{Dt}\o -(\p_1u\cdot\na) u_2+(\p_2u\cdot\na)u_1  = \frac{D}{Dt}\o + \o \div u , \ea
\ee
and that
\be\la{rb3}\ba \div  \dot u&=\frac{D}{Dt}\div u +(\p_1u\cdot\na) u_1+(\p_2u\cdot\na)u_2\\&=
 \frac{D}{Dt}\div u - 2\nabla u_1\cdot\nabla^{\perp}u_2+ (\div u)^2
\\&=
\frac{D}{Dt}(\frac{G}{2\mu+\lam})+ \frac{D}{Dt}(\frac{P-\bar P}{2\mu+\lam}) - 2\nabla u_1\cdot\nabla^{\perp}u_2 + (\div u)^2.\ea
\ee

Then, we rewrite
the momentum equations  as
\be\la{n5}
\rho\dot{u} = \nabla G + \mu\nabla^{\perp}\o.
\ee Multiplying  \eqref{n5}    by $2\dot u $  and integrating the resulting equality over $\O,$  we obtain after using \eqref{ra3} and \eqref{rb3} that
\be\la{p3}
\ba
& \frac{\rm d}{{\rm d}t} \int \left(\o^2 + \frac{G^2}{2\mu+\lam}\right)dx  + 2A_2^2\\
& = -\int \o^2\div udx + 4\int G\nabla u_1\cdot\nabla^{\perp}u_2dx- 2\int G(\div u)^2dx\\
 &\quad -\int\frac{(\beta-1)\lam - 2\mu}{(2\mu+\lam)^2}G^2\div udx + 2\beta \int\frac{\lam(P-\bar P)}{(2\mu+\lam)^2}G\div udx\\
 & \quad-2\g\int\frac{ P}{2\mu+\lam}G\div udx + 2(\g-1) \int P \div udx \int\frac{G}{2\mu+\lam}dx \\
 &  = \sum_{i=1}^7I_i.
 \ea
\ee

Each $I_i$ can be estimated as follows:

First,
it follows from \eqref{n5} that \be\la{a39} \lap G=\div(\n \dot u),\quad \mu\lap \o=\na^\perp\cdot (\n\dot u).\ee
which together with the standard $L^p$-estimate of elliptic equations implies that for $p\in (1,\infty),$\be \la{3.36}\|\na G\|_{L^p}+\|\na \o\|_{L^p}\le C(p,\mu)\|\n \dot u\|_{L^p}.\ee In particular, we have
\be\la{d3.19}  \|\na G\|_{L^2}+\|\na \o\|_{L^2}\le C(\mu)\rt^{1/2}A_2.\ee
This combining with \eqref{lp} gives
\be\la{a3.26}
\ba   \|\o\|_{L^4} & \le  C\|\o\|_{L^2}^{1/2}\|\na\o\|_{L^2}^{1/2}\\&\le   C \rt^{1/4}
 A_1^{1/2} A_2^{1/2},
\ea\ee
which leads to
\be\la{a3.27}\ba
|I_1| &\le C\|\o\|_{L^4}^2\|\div u\|_{L^2}  \\&\le \ve A_2^2+C(\ve)\rt A_3^{2}\vp^2,\ea
\ee for $\varphi$ as in \eqref{a3.15}.

Next, we will use an idea due to Perepelitsa \cite{Mik} to estimate $I_2.$ Noticing that $${\rm rot} \na u_1=0,\quad \div \nabla^{\perp}u_2=0,$$ one thus derives   from \cite[Theorem II.1]{coi3} that $$ \|\na u_1\cdot \nabla^{\perp}u_2\|_{\mathcal{H}^1}\le C\|\na u\|_{L^2}^2.$$
This combining with the fact that $\mathcal{BMO}$ is the dual space of $\mathcal{H}^1$ (see \cite{fef1}) gives
\be \la{a3.28}\ba |I_2|&\le C\|G\|_{\mathcal{BMO}}\|\na u_1\cdot \nabla^{\perp}u_2\|_{\mathcal{H}^1}\\&\le C\|\na G\|_{L^2}\|\na u\|_{L^2}^2\\&\le C \|G\|_{H^1}A_3\|\na u\|_{L^2} \\&\le C \|G\|_{H^1}A_3 \vp,  \ea\ee
where in the last inequality, we have used the following simple fact:
\bnn \la{a3.29}\ba  \|\na u\|_{L^2} &\le  C\|\o\|_{L^2}  +C\|\div u\|_{L^2} \\&\le C\|\o\|_{L^2} +C\left\|\frac{G}{2\mu+\lam}\right\|_{L^2} +C\left\|\frac{P-\bar P}{2\mu+\lam}\right\|_{L^2}  \le C \vp.\ea\enn

Next, Holder's inequality yields that
\be\la{pp4}\ba
 \sum\limits_{i=3}^6|I_i|
 &\le C\int|G||\frac{G+P-\bar P}{2\mu+\lam}||\div u|dx+C\int \frac{G^2}{2\mu+\lam} |\div u|dx\\&\quad +C\int \frac{P |G| }{2\mu+\lam} |\div u|dx+C\int \frac{ |G| |\div u| }{2\mu+\lam}dx\\
 &\le C\int \frac{G^2 |\div u|}{2\mu+\lam}dx +C\int \frac{P |G| }{2\mu+\lam} |\div u|dx+C\int \frac{ |G| |\div u| }{2\mu+\lam}dx \\& \le C A_3  \left\|\frac{G^2}{2\mu+\lam}\right\|_{L^2} +C A_3 \|G\|_{L^4} \left\|\frac{ P }{(2\mu+\lam)^{3/2}}\right\|_{L^4}+C\|G\|_{L^2}  A_3   .\ea
\ee
It   follows from    \eqref{an7} that
\be\la{ggaa} \ba \|G\|_{L^2}  \le C \rt^{\beta/2} A_1 ,\ea\ee  which together with
the Holder inequality and \eqref{lp} yields that for $0<\al<1,$
\be\la{p7}\ba
\left\|\frac{G^2}{\sqrt{2\mu+\lam}}\right\|_{L^2}&\le C\left\|\frac{G}{\sqrt{2\mu+\lam}}\right\|_{L^2}^{1-\al}
\|G\|_{ L^{2(1+\al)/\al}}^{1+\al}\\ &\le C(\al)A_1^{ 1-\al }\|G\|_{L^2}^{\al}
\|G\|_{H^1}   \\ &\le C(\al)A_1 \rt^{ \al\beta/2 } \|G\|_{H^1} \\ &\le C(\al)  \|G\|_{H^1}\vp .\ea
\ee
This combining with  \eqref{pp4}  and \eqref{ggaa} gives
\be \la{p4}
 \sum\limits_{i=3}^6|I_i|   \le C(\al)  \|G\|_{H^1} A_3 \vp .
\ee

Next, Holder's inequality leads to
\be\la{p6}\ba
|I_7|&\le C \left\| \frac{P }{(2\mu+\lambda)^{1/2}}\right\|_{L^2}  A_3 \|G\|_{L^2} \\ &\le C \|G\|_{H^1}A_3 \vp. \ea
\ee

Finally, noticing that $\bar G$ satisfies
\be \la{bg}|\bar G| \le C\|\n\|_{L^\beta}^{\beta/2}A_3,\ee
we deduce from   the Poincar\'e-Sobolev  inequality and \eqref{d3.19} that
\be\la{335}\ba   \|\o\|_{H^1}+\|G\|_{H^1}&\le C\|\na\o\|_{L^2}+C   \|G-\bar G\|_{L^2}+C|\bar G|+C\|\na G\|_{L^2}  \\& \le C\|\na\o\|_{L^2}+ C \|\na G\|_{L^2}+C\|\n\|_{L^\beta}^{\beta/2}A_3\\& \le C \rt^{1/2}A_2+C\|\n\|_{L^\beta}^{\beta/2}A_3.\ea\ee
Substituting \eqref{a3.27},  \eqref{a3.28},  \eqref{p4},    \eqref{p6},  and  \eqref{335} into  \eqref{p3}, we obtain that for any $\ve>0,$ \bnn
  \ba
   \frac{d}{dt}A_1^2(t) + 2A_2^2(t)  &\le \ve A_2^2+C(\ve)\rt A_3^{2} \vp^2 +C(\al)\|G\|_{H^1}A_3\vp \\&\le 2\ve A_2^2+C(\ve,\al)\rt A_3^{2} \vp^2 +C(\al)\|\n\|_{L^\beta}^{\beta/2}A_3^2\vp ,
  \ea
  \enn
 which directly gives  \eqref{a3.14} after choosing $\ve$ suitably small. The proof of Lemma \ref{le1} is completed.

Lemma \ref{le1} directly yields that
\begin{lemma}\la{co-1}
   For any $\al\in(0,1),$   there is a constant $C (\al) $ depending only on $\al,\mu,\beta,\gamma, $ $\|\n_0\|_{L^\infty}, $ and $\|u_0\|_{H^1} $   such that
\be\la{n9}
\ba
  &\sup\limits_{0\le t\le T}\log(e+A_1^2(t)+A_3^2(t)) + \int_0^T\frac{A_2^2(t)}{e + A_1^2(t)}dt    \le   C(\al)  \rt^{1+\ka+\al\beta},
\ea
\ee
 with
\be\la{342}\ka = \max\{0,\,\, (3\ga-6\beta) /2,\,\,\ga-\beta ,\,\, \beta-\ga-2\}.\ee
\end{lemma}

{\it Proof.} It follows from  \eqref{a3.15} and \eqref{r3} that
\be\la{var1} \ba
\vp(t)\le  C+C A_1 \rt^{ \al\beta /2 }    +  C\rt^{ (3\ga-6\beta) /4 } +  C\rt^{ ( \ga- \beta) /2 } , \ea
\ee which together with \eqref{a3.14}    and \eqref{r3}  gives \be\la{var31}
  \ba&
   \frac{d}{dt}A_1^2(t) +  A_2^2(t)    \\ &\le  C( \al)\rt A_3^{2} \vp^2 +C(\al)\|\n\|_{L^\beta}^\beta\rt^{-1} A_3^2 \\  &\le  C( \al)\rt\left(\rt^{\al\beta} A_1^2+  \rt^{(3\ga-6\beta) /2} +   \rt^{ \ga- \beta  }+\rt^{\beta-\ga-2} +1\right) A_3^2  .
  \ea
  \ee
   Dividing \eqref{var31} by $e+A_1^2(t)$ and integrating the resulting inequality over $(0,T), $  we obtain  \eqref{n9} after using  \eqref{c3.10}, \eqref{r3}, and  \eqref{a3.11}. We thus finish the proof of Lemma \ref{co-1}.

\begin{remark} Under the stringent condition  \eqref{pe1},   Perepelitsa \cite{Mik}  also obtained \eqref{n9} with $ \ka =0 .$
\end{remark}

  Next, we denote  the  commutator $F$   by \be \la{a3.42}F\triangleq \sum_{i,j=1}^2[u_i,R_iR_j](\rho u_j).\ee The following lemma   gives an   estimate of $F$ which will play an important role in obtaining the uniform upper bound of the density.
\begin{lemma}\la{lce} For any $\ve>0,$    there is a positive constant  $C (\ve)$ depending only on $\ve,  $ $     \mu, $ $   \beta,   $ $  \gamma,  $ $    \|\n_0\|_{L^\infty},$ and $ \| u_0\|_{H^1} $    such that
  \be\la{3.41}
  \ba
   \left\|F\right\|_{L^\infty}\le  \frac{ C (\ve)\rt^{-1-\ka} A^2_2}{ e+A_1^2 }+C (\ve) A_3^2 \rt^{(3+\ka)/2+\ve}  +C (\ve) \rt^{1+\ve} ,
  \ea
  \ee
  with $\ka$ as in \eqref{342}.

\end{lemma}

{\it Proof.} First, it follows from \eqref{r3} that \be\la{3.45}\ba \|\n u\|_{L^{2\ga/(\ga+1)}}&\le \|\n\|_{L^\ga}^{1/2}\|\n^{1/2}u\|_{L^2}\le C,\ea\ee which together with \eqref{u-1} gives
\be\la{3.46}\ba |\ol F|  \le C\|u\|_{L^{2\ga/(\ga-1)}}\|\n u\|_{L^{2\ga/(\ga+1)}} \le C+CA_3.\ea\ee

Then, we  deduce from \eqref{2.6} that
 \bnn \|F\|_{L^q}\le C(q)\|u\|_{\mathcal{BMO}}\|\n u\|_{L^q}\le C(q)\|\na u\|_{L^2}\|\n u\|_{L^q},\enn which together with
the Gagliardo-Nirenberg  inequality and  \eqref{2.7} thus gives that for $q\in (8,\infty),$
\be\la{a3.46} \ba\|F-\bar F\|_{L^\infty}&\le  C(q) \|F-\bar F\|_{L^q}^{ (q-4)/q}\|\na F\|_{L^{4q/(q+4)}}^{4/q}\\&\le   C(q)\left(\|\na u\|_{L^2} \|\n  u\|_{L^q} +|\bar F| \right)^{ (q-4)/q}\left(\|\na u\|_{L^4} \|\n  u\|_{L^q} \right)^{4/q}\\ &\le   C(q)A_3^{ (q-4)/q}\|\na u\|_{L^4}^{4/q} \|\n  u\|_{L^q} +C(q)\rt^{4/q} \left(A_3 +1\right)\|\na u\|_{L^4}^{4/q}  ,\ea \ee
where in the last inequality, we have used \eqref{3.46} and   the following simple fact:\bnn \|\n u\|_{L^q}\le C\rt \|u\|_{L^q}\le C(q)\rt (1+A_3)\enn due to \eqref{u-1}.

Next, noticing that \eqref{a3.11} gives
\bnn\ba e+A_1\le C +CA_3+C\rt^{(\ga-\beta)/2}\le C\rt^{\max\{0,(\ga-\beta)/2\}}(e+A_3),\ea\enn
we obtain from  \eqref{a3.26},  \eqref{p7},    \eqref{335}, and \eqref{var1}  that
\be\la{3.43}
\ba \|\na u\|_{L^4}   &
 \le C(\|\div u\|_4+\|\omega\|_4)\\ &\le  C\left\|\frac{G+P-\bar P}{
2\mu+\lambda}\right\|_{L^4} +  C \rt^{1/4}
 A_1^{1/2} A_2^{1/2}\\ &\le  C\left\|\frac{G^2 }{
\sqrt{2\mu+\lambda}}\right\|_{L^2}^{1/2} +C \left\|\frac{ P-\bar P}{
2\mu+\lambda}\right\|_{L^4} +  C \rt^{1/4}
 A_1^{1/2} A_2^{1/2}\\& \le C \left(\rt^{1/4}A_2^{1/2}+C\|\n\|_{L^\beta}^{\beta/4}A_3^{1/2}\right)\varphi^{1/2}
 +C\rt^{(3\ga-4\beta)/4}\\&\quad+  C \rt^{(5+\ka)/4 }  (e+A_1)  \left(\frac{\rt^{-4-\ka}A_2^2}{e+A_1^2}\right)^{1/4}\\
& \le C  \rt^{\ti C }  (e+A_3)  \left(\frac{\rt^{-4-\ka}A_2^2}{e+A_1^2}\right)^{1/4} + C  \rt^{\ti C } (e+A_3) ,\ea\ee with some constant $\ti C>1 $ depending only on $\beta$ and $\ga.$
This combining with \eqref{n9}  implies that for $\al \in(0,1),$
\bnn \la{3.44}\ba & \log
 \left(e+\|\na u\|_{L^4}\right)\\
& \le C(\al)  \log(e+ \rt )+ C \log(e+ A_3 )  + C \log(e+ \frac{\rt^{-4-\ka}A_2^2}{ (e+A_1^2) (e+A_3)^6} )\\
& \le  C(\al) \rt^{1+\ka+\al\beta}   +\frac{C\rt^{-4-\ka}A_2^2}{ (e+A_1^2) (e+A_3)^6}, \ea\enn
which together with    \eqref{bmo} and   \eqref{u-1} gives that for $\alpha \in (0,1),$
\be\la{x3.47}\ba  \|  u\|_{L^\infty}& \le C   \| \na u\|_{L^2} \log^{1/2
}\left(e+\|\na u\|_{L^4}\right)+C\|\na u\|_{L^2}+C \\&\le  C (\al) A_3 \rt^{(1+\ka+\al\beta)/2}+  C\left(\frac{\rt^{- 4-\ka }A_2^2}{ (e+A_1^2) (e+A_3)^4} \right)^{1/2} +C .\ea\ee
It thus follows from   Holder's inequality, \eqref{x3.47}, and \eqref{r3}   that for $\alpha \in (0,1)$ and  $q\in (8,\infty),$
\be\la{3.47}\ba \|\n u\|_{L^q}   &\le C \rt^{1-1/q
}  \|\n^{1/2}u\|_{L^2}^{2/q}\|u\|_{L^\infty}^{1-2/q
}\\ &\le C (\al)A_3^{1-2/q
}\rt^{ (3+\ka+\al\beta)/2
}   +C     \left(\frac{\rt^{-1-\ka}A_2^2}{ (e+A_1^2) (e+A_3)^4}\right)^{1/2-1/q
}+C\rt.\ea\ee

Substituting  \eqref{3.47} and \eqref{3.43} into  \eqref{a3.46} yields that  for $\alpha \in (0,1)$ and  $q\in (8,\infty),$
\be \la{3.48}\ba &\|F-\bar F\|_{L^\infty} \\ &\le C(q,\al) \|\na u\|_{L^4}^{4/q} A_3^{ (2q-6)/q}
  \rt^{(3+\ka+\al\beta)/2
} \\&\quad  + C(q,\al)A_3^{ (q-4)/q}\|\na u\|_{L^4}^{4/q}
   \left(     \frac{\rt^{-1-\ka}A_2^2}{ (e+A_1^2) (e+A_3)^4} \right)^{1/2-1/q
}  \\&\quad+C(q)\rt  A_3^{(q-4)/q}  \|\na u\|_{L^4}^{4/q}+C(q)\rt^{4/q} \left(A_3 +1\right)\|\na u\|_{L^4}^{4/q}\\ &
  \le C(q,\al)  A_3^{ (2q-6)/q}  (e+A_3)^{4/q} \left(\frac{\rt^{-1-\ka}A_2^2}{e+A_1^2}\right)^{1/q}\rt^{ (3+\ka+\al\beta)/2+4\ti C/q }\\&\quad +  C(q, \al)A_3^{ (2q-6)/q}  (e+A_3)^{4/q}  \rt^{ (3+\ka+\al\beta)/2+4\ti C/q }   \\& \quad+ C(q,\al) \|\na u\|_{L^4}^{4/q}
   \left(\frac{ \rt^{-1-\ka}A_2^2}{  e+A_1^2 }  \right)^{1/2-1/q
} +C(q)\rt  \left(A_3 +1\right)\|\na u\|_{L^4}^{4/q}\\&\triangleq \sum\limits_{i=1}^4J_i.  \ea\ee
Holder's inequality implies that
\be\la{3.49}\ba |J_1|& \le C(q,\al)\left(\rt^{ (3+\ka+\al\beta)/2+4\ti C/q }  A_3^{ (2q-6)/q}\right)^{q/(q-3)}  +(e+A_3^2) +\frac{\rt^{-1-\ka}A_2^2}{e+A_1^2}\\ &\le C(q,\al)+C(q,\al)  A_3^2 \rt^{ \ti\ka(\al,q) }  +\frac{\rt^{-1-\ka} A^2_2}{e+A_1^2} , \ea\ee
with
\be \la{3.50}\ti\ka(\al,q)\triangleq \left(\frac32+\frac{\ka }{2}+\frac{\al\beta }{2}+\frac{4\ti C}{q}\right)\frac{q}{q-3}. \ee
Similarly, we have
   \be \la{3.52}\ba |J_2|\le C(q,\al)+C(q,\al)  A_3^2 \rt^{ \ti\ka(\al,q) } .\ea\ee

One thus deduces from \eqref{3.43}  that for $\eta\in (0,1),$
\be\la{a3.53}
\ba
   \|\na u\|_{L^4}^{\eta} & \le     C(  \eta)   \rt^{4\ti C \eta/(4-\eta)}  (e+A_3)^{4  \eta/(4-\eta)}+    \frac{\rt^{-4-\ka}A_2^2}{e+A_1^2}   \\ &\quad+ C(\al,\eta) \rt^{\ti C \eta} ( A_3^2 +1 )\\ &\le  C( \eta) \rt^{2\ti C \eta} ( A_3^2 +1 )+    \frac{\rt^{-4-\ka}A_2^2}{e+A_1^2}  ,\ea\ee
which together with Holder's inequality gives
\be\la{3.55}
\ba |J_3| &\le  \|\na u\|^{8/q}_{L^4}+  \frac{\rt^{-1-\ka}A_2^2}{e+A_1^2}   +C(q,\al) \\&\le       C(\al,q) \rt^{16\ti C/q} + C(\al,q)\rt^{16\ti C/q}  A_3^{2} +  \frac{ 2\rt^{-1-\ka} A^2_2}{e+A_1^2 }.\ea\ee

It follows from \eqref{u-1}  and \eqref{a3.53} that\be\la{3.53}\ba |J_4|  &\le  \rt \|\na u\|_{L^4}^{8/q} +C(q)\rt +C(q)\rt  A_3^2 \\&\le       C( q) \rt^{1+16\ti C/q} + C( q)\rt^{1+16\ti C/q}  A_3^{2} +  \frac{ \rt^{-1-\ka} A^2_2}{e+A_1^2 }.\ea\ee

Substituting \eqref{3.49}, \eqref{3.52}, \eqref{3.55},  and \eqref{3.53} into \eqref{3.48},   we obtain \eqref{3.41} after using \eqref{3.46} and choosing $q$ suitably large and then $\al$ suitably small.
  The proof of Lemma \ref{lce} is completed.

Now, we are in a position to prove Proposition \ref{lle5}.

 {\it Proof of Proposition \ref{lle5}. }
It follows from \eqref{n5} that  $G$ solves
\bnn\la{n6}
\lap G = \div(\rho\dot{u})=\p_t(\div(\rho u)) + \div\div(\rho u\otimes u),
\enn
which implies
\be\la{n12}
\ba
  &G -\bar G+ \frac{D}{Dt}\left((-\lap)^{-1}\div(\rho u)\right)=F,
\ea
\ee with $F$ as in \eqref{a3.42}.
 The mass equation $(\ref{n1})_1$ leads to \bnn -\div u =\frac{1}{\n}D_t\n ,\enn
  which
  combining with (\ref{n12})  gives that
\be\la{key}
 \frac{D}{Dt}\theta(\rho)  + P = \frac{D}{Dt}\psi + \bar P-\bar G +F ,
\ee
with \be\la{p26}
\theta(\n) \triangleq 2\mu\log\rho +  \beta^{-1}\rho^{\beta}  ,\quad
 \psi\triangleq(-\lap)^{-1} \div(\rho u).\ee
Since the function $y=\theta(\n)$ is increasing for $\n\in (0,\infty),$ the inverse function \be \la{368}\n= \theta^{-1} (y)\ee exists for $y\in (-\infty, \infty).$  We rewrite \eqref{key} as \bnn \frac{D}{Dt}y = g(y)+\frac{D}{Dt}h,\enn
with \be\la{369} y=\theta(\n),\quad g(y)=-P(\theta^{-1}(y)),\quad h=\psi+\int_0^t\left( \bar P-\bar G+F\right)dt.\ee
 To apply Lemma \ref{ale1}, noticing that \bnn \lim\limits_{y\rightarrow \infty}g(y)=-\infty,\enn we need to  estimate $h.$ First, it follows from \eqref{bmo},    \eqref{3.45},   and \eqref{u-1} that
\be\la{psi}
\ba
\|\psi\|_{L^\infty} & \le C  \|{\na\psi}\|_{L^2} \log^{1/2}(e+ \|\nabla\psi \|_{L^3})+C\|\psi\|_{L^2}+C\\ & \le C  \| \n u\|_{L^2} \log^{1/2}(e+ \|\n u \|_{L^3})+C\|\n u\|_{L^{2\ga/(\ga+1)}}+C\\ & \le C \rt^{1/2}  \log^{1/2}(e+ \rt(1+\|\na u \|_{L^2}))+C\\ & \le C \rt^{1/2}\log^{1/2}(e+ A_3^2) +C\rt ,\ea\ee which together with \eqref{n9} gives that for $\ka$ as in \eqref{342} \be\la{apsi}\ba   \|\psi\|_{L^\infty}\le C \sup\limits_{0\le t\le T}\rt^{(3+\ka)/2}.
\ea
\ee

 Next,  on one hand, \eqref{r3} and \eqref{bg} lead to
\be\la{a372}\ba
|\bar P-\bar G | &\le C + C \|\n\|_{L^\beta}^{\beta/2}A_3(t) \\ &\le C +C A_3^2(t)+C A_3^2(t)\rt^{\beta-\ga}.\ea
\ee On the other hand, one deduces    from   \eqref{3.41},   \eqref{n9}, and \eqref{r3} that  for any $\ve>0$   and  all $0\le t_1\le t_2\le T$
 \be\la{p22}\ba \int_{t_1}^{t_2}\|F\|_{L^\infty}ds\le C(\ve)  \rt^{ (3+\ka)/2 +\ve}+C(\ve)  \rt^{1+\ve}  ({t_2}-{t_1}).\ea\ee
This combining with \eqref{apsi} and \eqref{a372} implies that for all $0\le t_1\le t_2\le T$ and any $\ve>0 ,$
\bnn\ba |h(t_2)-h(t_1)|\le& C(\ve )\rt^{\max\{(3+\ka)/2+\ve, \beta -\ga  \}} +C(\ve )\rt^{1+\ve} (t_2-t_1).\ea\enn
 Therefore, one can
choose $N_0$ and $N_1$ in (\ref{a100}) as:
\bnn N_0 = C(\ve)\rt^{\max\{(3+\ka)/2+\ve, \beta -\ga  \} },\quad N_1= C(\ve)\rt^{1+\ve}  .\enn
For $g(y)$ as in \eqref{369} with $\n=\te^{-1}(y)$ as in \eqref{368} being the inverse function of $y=\te(\n),$ we have
$$ g(\zeta)= - (\theta^{-1}(\zeta))^\ga\le -N_1
=-C(\ve)\rt^{1+\ve}
 ,$$   for all $\zeta\ge\ti\zeta\triangleq  C(\ve) \rt^{\beta(1+\ve) /\ga}.$ Lemma \ref{ale1} thus yields that
 \bnn    \rt^{\beta}\le C(\ve)\rt^{\max\{(3+\ka)/2+\ve, \beta -\ga,\beta(1+\ve)/\ga  \} }, \enn
which together with \eqref{1.11} gives \eqref{a3.56}. We finish the proof of Proposition \ref{lle5}.

The following Proposition \ref{caq1}, which will play  an important role in obtaining the large-time behavior of $(\n,u),$ is a direct consequence of \eqref{var31},   \eqref{a3.56}, \eqref{r3},   \eqref{u-1}, \eqref{335}, and Gronwall's inequality.

\begin{proposition}  \la{caq1}  If  \eqref{1.11} holds, there is a positive constant  $C$ depending only on $   \mu,  \beta, \gamma, \|\n_0\|_{L^\infty},$ and $\|u_0\|_{H^1} $ such that
\be\la{ab12}
\sup_{0\le t\le T}\left(\|\n\|_{L^\infty}+\norm[H^1]{ u}\right) + \int_0^T\left(\|\o\|_{H^1}^2+\|G\|_{H^1}^2+A_2^2(t) +A_3^2(t) \right)dt\le C.
\ee \end{proposition}

\subsection{Time-dependent upper bound of the density}

The following Proposition \ref{aupper}  will give a time-dependent   upper bound of the density which is the key to obtain higher order estimates provided \eqref{bet} holds. Throughout this subsection,   $C$ denotes a
generic positive constant
  independent of  the lower bound of the initial  density.
\begin{proposition}\la{aupper} Assume that \eqref{bet} holds. Then   there is a positive constant  $C (T)$ depending only on $T,  \mu,  \beta,$ $ \gamma,     \|\n_0\|_{L^\infty},$ and $\| u_0\|_{H^1} $ such that
   \be\la{b3.56}
\sup_{0\le t\le T}\left(\|\n\|_{L^\infty}+\norm[H^1]{ u}\right)  + \int_0^T\left(\|\o\|_{H^1}^2+\|G\|_{H^1}^2+A_2^2(t)\right)dt\le C(T).
\ee
\end{proposition}

 Before proving Proposition \ref{aupper}, we establish some a priori estimates, Lemmas \ref{jle2} and \ref{lemma1}.

 We first state the $L^p$-estimate of the density   due to Vaigant-Kazhikhov (\cite{Ka}).
\begin{lemma} [\cite{Ka}] \la{kq1} Let $\beta>1.$   For any $1\le p<\infty,$ there is a positive constant  $C(T) $ depending only on $  T,    \mu, \beta,  \gamma,   \|\n_0\|_{L^\infty}, $ and $\| u_0\|_{H^1} $  such that
  \be\la{ja1}
  \sup_{0\le t\le T}\norm[L^p]{\rho(\cdot, t)}\le C (T) p^{\frac{2}{\beta-1}}.
  \ee\end{lemma}

 The following  $L^p$-estimate of the momentum which plays an important role in the estimate of the upper bound of the density is a direct consequence of  Lemma \ref{kq1}.
\begin{lemma}  \la{jle2}  Let $\beta>1.$  For any $q>4,$  there is a positive constant  $C(q,T)$ depending only on $ T,q,   \mu, $ $  \beta, \gamma,   \|\n_0\|_{L^\infty},$ and $\|u_0\|_{H^1} $ such that
\be\la{jle3}\ba \|\n u\|_{L^q}   \le  C(q,T)\rt^{1+\beta(q-2)/(4q)}   (1+ A_3)^{ 1-2/q}.\ea\ee
  \end{lemma}

{\it Proof.} First, we claim that there is a  positive constant  $\nu_0\le 1/2$  depending only on $\mu$ such that
  \be\la{ja2}
  \sup_{0\le t\le T}\int \n |u|^{2+\nu}dx\le C (T) ,
  \ee
  with \bnn\la{ja3}\nu= \rt^{-\beta/2}\nu_0\in (0,1/2].\enn

  Then, let
$ r\triangleq (q-2)(2+\nu)/\nu>2$ due to $q>4.$ It follows from Holder's inequality, \eqref{ja2},  \eqref{lp}, and \eqref{u-1} that
\bnn\ba \|\n u\|_{L^q}  &\le C \|\n u\|_{L^{2+\nu}}^{2/q}\|\n u\|_{L^r}^{1-2/q
}\\ &\le C(T)\rt^{1-2/((2+\nu)q)
}   \|  u\|_{L^r}^{1-2/q
}\\ &\le  C(T)\rt  \left(  r^{1/2} \|u\|_{H^1}\right)^{1-2/q
}\\ &\le  C(q,T)\rt^{1 +\beta(q-2)/(4q)}  (1+ \|\na  u\|_{L^2})^{ 1-2/q},\ea\enn  which together with \eqref{a3.12} shows \eqref{jle3}.

Finally, it remains to prove \eqref{ja2}.  Multiplying $(\ref{n1})_2$ by $(2+\nu)|u|^\nu u,$  we get after  integrating the resulting equation over $ \O$ that
\bnn\ba & \frac{d}{dt}\int \n |u|^{2+\nu}dx+ (2+\nu) \int|u|^\nu \left(\mu |\na u|^2+(\mu+\n^\beta) (\div u)^2\right) dx \\& \le  (2+\nu)\nu  \int (\mu+\n^\beta)|\div u||u|^\nu |\na u|dx +C  \int \n^\ga |u|^\nu |\na u|dx\\& \le \frac{2+\nu}{2}\int (\mu+\n^\beta)(\div u)^2|u|^\nu  dx +\frac{2+\nu}{2}\nu_0^2 (\mu+1) \int  |u|^\nu |\na u|^2 dx \\&\quad+ \mu  \int   |u|^\nu |\na u|^2dx + C \int \n |u|^{2+\nu} dx+C\int \n^{(2+\nu)\ga-\nu/2}dx ,\ea\enn
which, after choosing $\nu_0(\mu)$ suitably small, together with Gronwall's inequality and \eqref{ja1} thus gives \eqref{ja2}. The proof of  Lemma \ref{jle2} is completed.

The next lemma will deal with the   time-dependent   estimate on the spatial $L^\infty$-norm of the commutator operator $F$ defined by \eqref{a3.42}.

\begin{lemma} \la{lemma1} Let $\beta>1.$   For any $\ve>0,$    there is a positive constant  $C (\ve,T) $ depending only on $ \ve, T,  \mu,  \beta, \gamma,$ $    \|\n_0\|_{L^\infty},$ and $\| u_0\|_{H^1} $ such that
  \be\la{a3.41}
  \ba
   \left\|F\right\|_{L^\infty}\le  \frac{ C (\ve,T) A^2_2}{e+A_1^2}+C (\ve,T)  (1+A_3^2 )\rt^{1+\beta/4+\ve}  .
  \ea
  \ee\end{lemma}

{\it Proof.} First, it follows from  \eqref{3.43} and   \eqref{jle3}  that for  $q>8,$
  \be\la{a3.63}\ba &A_3^{ (q-4)/q}\|\na u\|_{L^4}^{4/q} \|\n  u\|_{L^q}\\ & \le C(q,T)\rt^{1 +\beta(q-2)/(4q)}   \left(A_3^{ 2-6/q}+1\right) \|\na u\|_{L^4}^{4/q}   \\ & \le C(q, T)\rt^{1 +\beta(q-2)/(4q) +4\ti C/q }\left(A_3^{ 2-6/q}+1\right) (1+A_3^2)^{ 2/q} \left(\frac{ A_2^2}{e+A_1^2}\right)^{1/q}  \\&\quad + C(q  ,T)\rt^{1 +\beta(q-2)/(4q) +4\ti C/q } \left(A_3^{ 2-6/q}+1\right) (1+A_3^2)^{ 2/q} \\ & \le C(q,  T)\rt^{1 +\beta(q-2)/(4q) +4\ti C/q }\left(A_3^{ 2-2/q}+1\right) \left(\frac{ A_2^2}{e+A_1^2}\right)^{1/q}  \\&\quad + C(q, T)\rt^{1 +\beta(q-2)/(4q) +4\ti C/q } \left(A_3^2+1\right)  \\ & \le C(q, T)\rt^{1 +\beta/4 +4\ti C/(q-1) }  (1+A_3^2) + \frac{ A_2^2}{e+A_1^2}.  \ea\ee
This combining with \eqref{a3.46}  and \eqref{3.53} yields that
\bnn  \ba  \|F-\bar F\|_{L^\infty} & \le C(q ,T)\rt^{1 +\beta/4 +4\ti C/(q-1) }   A_3^2 + C(q ,T)\frac{ A_2^2}{e+A_1^2}\\& \quad+ C( q) \rt^{16\ti C/q} + C( q)\rt^{16\ti C/q}  A_3^{2}, \ea \enn which together with \eqref{3.46} directly gives \eqref{a3.41} after choosing $q$ suitably large and then $\al$ suitably small. The proof of Lemma \ref{lemma1} is completed.

{\it Proof of Proposition \ref{aupper}.}
We deduce from  \eqref{a3.15} and \eqref{ja1} that  for any $\al\in (0,1),$\bnn\la{var2} \ba
\vp(t)\le  C(T,\al)+C(T,\al) A_1 \rt^{ \al\beta /2 } ,\ea\enn which together with Lemmas \ref{le1} and   \ref{kq1} gives  \be\la{nn9}
\ba
& \sup\limits_{0\le t\le T}\log(e+A_1^2(t)) + \int_0^T\frac{A_2^2(t)}{e + A_1^2(t)}dt   \le    C(T,\al)  \rt^{1+\al\beta }  .
\ea
\ee

Then, for $\psi$   as in \eqref{p26}, it follows from \eqref{psi} and \eqref{nn9} that \bnn \|\psi\|_{L^\infty}\le C(T)\rt^{4/3},\enn
which together with \eqref{key}, \eqref{a372}, \eqref{ja1}, \eqref{a3.41},  and \eqref{nn9}  yields that  for $\ve\in(0,1),$
\bnn    \rt^{\beta}\le  C(\ve,T)\rt^{\max\{ 1+\beta/4+\ve,\,\beta-\ga,\,4/3\}}. \enn
Due to \eqref{bet}, after choosing $\ve$ suitably small, this directly gives \bnn
   \sup\limits_{0\le t\le T}\|\n\|_{L^\infty}\le C(T),
  \enn which together with \eqref{var31},    \eqref{r3}, \eqref{u-1}, \eqref{335}, and Gronwall's inequality yields  \eqref{b3.56}.  We complete the proof of Proposition \ref{aupper}.

\section{\la{se4} A priori estimates (II): higher order estimates}

\begin{lemma}\la{310} Assume that   \be \la{b3.57}\sup\limits_{0\le t\le T}\|\n\|_{L^\infty}\le M,\ee  for some positive constant $M.$ Then    there is a positive constant  $C (M)$ depending only on $ M,  \mu,  \beta, \gamma,   $ $ \|\n_0 \|_{L^\infty}, $ and $\|u_0\|_{H^1}  $ such that
\be\la{b19}
\sup_{0\le t\le T}\sigma\int \rho|\dot{u}|^2dx + \int_0^T\sigma\|\na\dot{u}\|_{L^2}^2dt\le C(M),
\ee with $\sigma
(t)\triangleq\min\{1,t\}.$ Moreover,  if \eqref{1.11} holds, for any $p\in [1,\infty),$ there is a positive constant  $C (p)$ depending only on $p,  \mu,  \beta, \gamma,   $ $ \|\n_0 \|_{L^\infty}, $ and $\|u_0\|_{H^1}  $ such that for any $T\in (1,\infty),$ \be\la{ab19}\sup_{1\le t\le T} \|\na u\|_{L^p} \le C(p) .
\ee  
\end{lemma}

{\it Proof.} Operating $\dot u^j[\pa/\pa t+\div
(u\cdot)]$ to $ (\ref{n1})_2^j,$ summing with respect
 to $j,$ and integrating the resulting equation over ${\O}$, one obtains
after integration by parts that
\be\la{b13}
 \ba &
\left(\frac{1}{2}\int\rho|\dot{u}|^2dx \right)_t\\
& =   -\int \dot{u}_j[\p_jP_t +\text{div}(\p_jPu)]dx +
\mu\int\dot{u}_j[\p_t\triangle u_j +
\text{div}(u\triangle u_j)] dx\\
&\quad  + \int\dot{u}_j[\p_{jt}((\mu+\lam)divu) +
 div(u\p_j((\mu+\lam)divu))]dx
 \triangleq\sum_{i=1}^{3}N_i. \ea \ee

  First, using the equation $(\ref{n1})_1,$ we obtain after integration
by parts  that
\be\la{b14}
 \ba
N_1 & = - \int \dot{u}_j[\p_jP_t + \text{div}(\p_jPu)]dx \\
& =  \int [-P^{'}\rho\text{div}u\p_j\dot{u}_j +
\p_k(\p_j\dot{u}_ju_k)P - P\p_j(\p_k\dot{u}_ju_k)]dx
\\
&\le C(M) \|\nabla u\|_{L^2}  \|\nabla \dot u\|_{L^2}\\
&\le \frac{\mu}{8} \|\nabla \dot u\|^2_{L^2} +C(M) \|\nabla  u\|^2_{L^2}  . \ea \ee Then, integration by parts leads to \be\la{ab14} \ba
N_2 & =  \mu\int\dot{u}_j[\p_t\triangle u_j
+ \text{div}(u\triangle u_j)]dx \\
& = - \mu\int \left(|\nabla\dot{u}|^2 + \p_i\dot{u}_j\p_ku_k\p_iu_j - \p_i\dot{u}_j\p_iu_k\p_ku_j - \p_iu_j\p_iu_k\p_k\dot{u}_j\right)dx \\
&\le -\frac{ 3\mu}{4} \int|\nabla\dot{u}|^2dx  + C(M) \int
|\nabla u|^4dx . \ea \ee
Similarly,
\be\la{b15}
\ba
 N_3 & = \int\dot{u}_j[\p_{jt}((\mu+\lam)\div u) +
 div(u\p_j((\mu+\lam)\div u))]dx\\
 & = -\int \p_j\dot{u}_j[ ((\mu+\lam)\div u)_t    +
 \div(u(\mu+\lam)\div u)] dx\\
 &\quad -\int \dot{u}_j\div(\p_ju(\mu+\lam)\div u) dx\\
 & \le -\int\p_j\dot{u}_j[(\mu+\lam)\div u_t+ \lam_t\div u  +
 (u\cdot\nabla\lam )\div u + (\mu+\lam)(u\cdot\nabla) \div u]dx \\
  &\quad+  \frac{\mu}{8}\int|\nabla\dot{u}|^2dx +
 C(M)\int|\nabla u|^4dx\\
&  = -\int\left(\frac{D}{Dt}\div u+ \pa_ju_i\pa_iu_j\right)[(\mu+\lam)\frac{D}{Dt}\div u -\n\lam'(\n) \div u]dx\\ & \quad+  \frac{\mu}{8}\int|\nabla\dot{u}|^2dx
  +C(M) \int|\nabla u|^4dx\\
& \le -\frac{\mu}{2}\int|\frac{D}{Dt}\div u|^2dx   +  \frac{\mu}{8}\int|\nabla\dot{u}|^2dx + C(M)\|\nabla   u\|^4_{L^4}+C(M)\|\nabla   u\|^2_{L^2} .
\ea
\ee
Finally, substituting \eqref{b14}-\eqref{b15} into  \eqref{b13}  shows that
\be\la{b17}
 \ba & 2\left(
\int\rho|\dot{u}|^2dx \right)_t + \mu\int
|\nabla\dot{u}|^2dx+\mu \int |\frac{D}{Dt}\div u|^2dx\\
&\le C(M)|\na u\|_{L^4}^4 +C(M)\|\nabla   u\|^2_{L^2}   \\ &\le C(M)(\norm[L^4]{G}^4 + \norm[L^4]{\o}^4 + \norm[L^4]{P-\ol P}^4+\|\nabla   u\|^2_{L^2})\\
 &\le C(M)(\norm[L^2]{G}^2\norm[H^1]{ G}^2 + \norm[L^2]{\o}^2\norm[L^2]{\na\o}^2 + \norm[L^2]{P-\ol P}^2+\|\nabla   u\|^2_{L^2})\\
 &\le C(M)( \norm[H^1]{ G}^2 +  \norm[L^2]{\na\o}^2  +\|\nabla   u\|^2_{L^2})\\
 & \le C(M)\norm[L^2]{\rho^{\frac{1}{2}}\dot{u}}^2 + C(M)\|\nabla   u\|^2_{L^2},
\ea
\ee where in the last inequality we have used \eqref{335} and \eqref{a3.12}.
Multiplying \eqref{b17} by  $\sigma,$ integrating the resulting equation over $(0,T),$ we obtain (\ref{b19}) after using \eqref{b3.56}.

It remains to prove \eqref{ab19}. Because of \eqref{1.11}, we deduce from \eqref{1.11} and \eqref{a3.56}  that \bnn\ba \|\na u\|_{L^p} &\le C(p)  \|\div u\|_{L^p}+C(p)  \|\o \|_{L^p}\\& \le  C(p)  \|G\|_{L^p}+ C(p)  \|P-\ol P\|_{L^p}+C(p)  \|\o \|_{L^p}\\& \le  C(p)  \|G\|_{H^1} +C(p)  \|\o\|_{H^1}+ C(p)\\& \le  C(p)  \|\n^{1/2}\dot u\|_{L^2} +C(p) ,   \ea\enn where in the last inequality we have used \eqref{335}, \eqref{a3.12}, and \eqref{ab12}. This combining with \eqref{b19} gives \eqref{ab19}.  We finish the proof of Lemma \ref{310}.

\begin{lemma}\la{d31} Assume that \eqref{bet} holds. Then  for any $p>2,$  there is a positive constant  $C (p,T )$ depending only on $ p, T,  \mu,  \beta, \gamma,   $ $ \|\n_0 \|_{L^\infty}, $ and $\|u_0\|_{H^1}  $ such that
\be \la{4a2} \ba&\int_0^T \left(\|G\|_{L^\infty}+\|\na G\|_{L^p}+\|\na \o\|_{L^p}+\| \rho \dot u\|_{L^p}\right)^{1+1 /p}dt \\& +\int_0^Tt\left(\|\na G\|_{L^p}^2+\|\na \o\|_{L^p}^2+\|   \dot u\|_{H^1}^2\right) dt \le C(p,T) .\ea\ee
\end{lemma}

{\it Proof. }
If follows from \eqref{2.5}, \eqref{lp}, and \eqref{b3.56} that
 \bnn\la{b22}\ba
 \| \rho \dot u\|_{L^p} & \le
 C\| \rho \dot u\|_{L^2}^{2(p-1)/(p^2-2)}\|\dot u\|_{L^{p^2}}^{p(p-2)/(p^2-2)}\\ & \le
 C\| \rho \dot u\|_{L^2}^{2(p-1)/(p^2-2)}\|\dot u\|_{H^1}^{p(p-2)/(p^2-2)}\\ & \le
 C\| \rho  \dot u\|_{L^2} +C\| \rho \dot u\|_{L^2}^{2(p-1)/(p^2-2)}\|\na \dot u\|_{L^2}^{p(p-2)/(p^2-2)}   , \ea\enn which together with   \eqref{b3.56},  \eqref{b19},   and \eqref{2.5} implies that\be\la{gh1}  \ba&\int_0^T \left(\| \rho \dot u\|^{1+1 /p}_{L^p}+t\|  \dot u\|^2_{H^1}\right) dt\\ &\le C(p,T) \int_0^T\left( \| \rho^{1/2}  \dot u\|_{L^2}^2 +  t\|\na \dot u\|_{L^2}^2+ t^{-(p^3-p^2-2p)/(p^3-p^2-2p+2)} \right)dt \\ &\le C(p,T) .\ea\ee
Noticing  that the  Gargliardo-Nirenberg inequality and  \eqref{b3.56} yield that
 \be\la{419}\ba &\|\div u\|_{L^\infty}+\|\o\|_{L^\infty} +  \|G\|_{L^\infty} +\|\o\|_{L^\infty}\\ &\le C(p,T) +C(p,T) \|\na G\|_{L^p}^{p/(2(p-1))} +C(p,T) \|\na \o\|_{L^p}^{p/(2(p-1))}\\ &\le C(p,T) +C(p,T) \|\n\dot u\|_{L^p}^{p/(2(p-1))} , \ea\ee we directly derive \eqref{4a2}  from \eqref{gh1} and \eqref{3.36}.  We finish the proof of Lemma \ref{d31}.

\begin{proposition}\la{le5} Assume that \eqref{bet} holds.   Then,  for   $q>2$  as in Theorem \ref{t2}, there is a constant $C(T)$ depending only on  $ T,
q, \mu,  \ga,\beta, \| u_0\|_{H^1},$ and $
    \|\n_0 \|_{W^{1,q}}$ such that
  \be\la{pa1}\ba
  &\sup_{0\le t\le T}\left(\norm[W^{1,q}]{ \rho}+\| u\|_{H^1}+  t\| u\|_{H^2}  \right)\\&+\int_0^T \left(\|\nabla^2 u\|_{L^q}^{(q+1)/q}+t\|\nabla^2 u\|_{L^q}^2+t \|u_t\|_{H^1}^2\right)dt\le C(T).\ea
  \ee
\end{proposition}

{\it Proof.}  Following \cite{hlx}, we will   prove \eqref{pa1}.
First,
 denoting by  $\Phi\triangleq  (\Phi^1,\Phi^2)$ with $\Phi^i\triangleq (2\mu+\lambda(\n))\pa_i \n\,\, ( i=1,2),$ one deduces from $(\ref{n1})_1$ that $\Phi^i$ satisfies
\be \la{4.20}\ba  &\Phi^i_t +(u\cdot\na)  \Phi^i +(2\mu+\lambda(\n)) \pa_i u^j \pa_j \n +  \n\pa_i G+  \n\pa_i   P  +  \Phi^i \div u  =0 .\ea\ee
For $ q> 2 ,$ multiplying \eqref{4.20} by $|\Phi |^{q-2}  \Phi^i $ and integrating the resulting equation over $\O,$ we obtain after integration by parts and using \eqref{3.36} that
  \be\la{b21}\ba
\frac{d}{dt}\norm[L^q]{\Phi}   & \le
 C(1+\norm[L^{\infty}]{\nabla u} )
\norm[L^q]{\nabla\rho} +C\| \na G\|_{L^q} \\ & \le
 C(1+\norm[L^{\infty}]{\nabla u} )
\norm[L^q]{\nabla\rho} +C\| \rho \dot u\|_{L^q}. \ea\ee

  Next,
   we deduce from  standard $L^p$-estimate for elliptic  system,   \eqref{419}, and  \eqref{3.36}  that
   \be\la{420}\ba \|\na^2u\|_{L^q}&\le C\|\na\div u\|_{L^q}+C\|\na \o\|_{L^q}\\ &\le C\|\na((2\mu+\lambda)\div u)\|_{L^q}+C \|\div u\|_{L^\infty} \|\na \n\|_{L^q}+C\|\na \o\|_{L^q}\\ &\le   C(\|\div u\|_{L^\infty}+1)\|\na \n\|_{L^q}+C\|\na G\|_{L^q}+C\|\na \o\|_{L^q}\\ &\le   C(\|\n\dot u\|_{L^q}^{q/(2(q-1))}+1)\|\na \n\|_{L^q}+C\|\n\dot u\|_{L^q} \\ &\le  C\|\na \n\|_{L^q}^{(2q-2)/(q-2)}+C\|\n\dot u\|_{L^q} +C .\ea\ee
   
  Then, it follows from
Lemma \ref{le9}, \eqref{419}, and \eqref{420} that
    \be\la{b24}\ba   \|\na
u\|_{L^\infty }  &\le C\left(\|{\rm div}u\|_{L^\infty }+
\|\o\|_{L^\infty } \right)\log(e+\|\na^2 u\|_{L^q}) +C\|\na
u\|_{L^2} +C \\&\le C\left(1+\|\n\dot u\|_{L^q}^{q/(2(q-1))}\right)\log(e+\|\rho \dot u\|_{L^q} +\|\na \rho\|_{L^q}) +C\\&\le C\left(1+\|\n\dot u\|_{L^q} \right)\log(e+   \|\na \rho\|_{L^q}) . \ea\ee

 Substituting \eqref{b24} into \eqref{b21}, we deduce from Gronwall's inequality and \eqref{4a2}  that \bn \la{b30} \sup\limits_{0\le t\le T}\|\nabla
\rho\|_{L^q}\le C,\en  which combining with  \eqref{420} and \eqref{4a2} shows
  \be \la{b31}\int_0^T \left(\|\nabla^2 u\|_{L^q}^{(q+1)/q}+t\|\nabla^2 u\|_{L^q}^2\right)dt\le C.\ee
Finally, it  follows from  \eqref{2.5}, \eqref{b3.56},  \eqref{b19}, and \eqref{b31} that
\be\la{po1}\ba &\int_0^Tt \|u_t\|_{H^1}^2dt\\ &\le C\int_0^Tt \left(\| \n^{1/2} u_t\|_{L^2}^2+  \|\na u_t\|_{L^2}^2\right)dt\\&\le C\int_0^Tt \left(\|\n^{1/2}\dot u\|_{L^2}^2 +  \| u\cdot\na u  \|_{L^2}^2+\|\na\dot u \|_{L^2}^2 +  \|\na(u\cdot\na u) \|_{L^2}^2\right)dt\\ &\le C +C\int_0^Tt \|\na u\|_{L^4}^4dt+C\int_0^Tt \|  u\|_{H^1}^2\|\na^2 u\|_{L^q}^2dt\\ &\le C +C\int_0^Tt \|\na u\|_{L^2}^2 \|\na^2 u\|_{L^2}^2dt+C\int_0^Tt  \|\na^2 u\|_{L^q}^2dt\\ &\le C.\ea\ee
   We obtain from   \eqref{b3.56},  \eqref{3.36}, and  \eqref{b30}  that
   \be \la{4.19}\ba \|\na^2u\|_{L^2}&\le C\|\na \o\|_{L^2} +C\|\na\div u\|_{L^2}\\ &\le C\|\na \o\|_{L^2} +C\|\na((2\mu+\lambda)\div u)\|_{L^2}+C \|\div u\|_{L^{2q/(q-2)}} \|\na \n\|_{L^q} \\ &\le  C\|\na \o\|_{L^2} +C\|\na G\|_{L^2}+C+ C \|\na u\|_{L^2}^{(q-2)/q}\|\na^2 u\|_{L^2}^{2/q} \\ &\le   C+\frac12  \|\na^2 u\|_{L^2} +C\|\n \dot u\|_{L^2} ,\ea\ee
which together with \eqref{b19} gives \bnn\sup\limits_{0\le t\le T}t \|\na^2u\|_{L^2}\le C. \enn
 This combining with
  \eqref{b30}--\eqref{po1} and \eqref{b3.56}  yields   \eqref{pa1}.
The proof of Proposition \ref{le5} is completed.

\section{ Proofs of  Theorems  \ref{t2}--\ref{th2}}

With all the a priori estimates in Sections \ref{se3} and \ref{se4}
at hand, we are ready to prove the main results of this paper in
this section. We first state the global existence of strong solution $(\n,u)$ provided that \eqref{bet} holds and  that $(\rho_0,m_0
)$ satisfies  \eqref{2.1}.

\begin{proposition} \la{pro2}
Assume that  \eqref{bet} holds and that  $(\rho_0,m_0
)$ satisfies  \eqref{2.1}.  Then   there exists a unique strong solution  $(\rho,u ) $      to   \eqref{n1}-\eqref{n4}
 in $\O\times (0,\infty)$     satisfying  \eqref{2.2}  for  any $T\in (0,\infty).$  In addition,    $(\rho,u ) $  satisfies \eqref{pa1} with some positive constant $C$ depending only on $T,$  $  \mu,  \beta, \gamma, \|\n_0\|_{L^\infty},$ and $\|u_0\|_{H^1}  $   such that.    Moreover, if \eqref{1.11} holds, there exists some positive constant $C$ depending only on $  \mu,  \beta, \gamma, \|\n_0\|_{L^\infty},$ and $\|u_0\|_{H^1} $ such that both \eqref{ab12} and  \eqref{ab19}  hold.

 \end{proposition}

{\it Proof.}
First,  standard local existence result, Lemma \ref{th0},   applies to show that the  problem   \eqref{n1}-\eqref{n4} with   initial data $(\n_0 ,m_0 )$ has   a unique local solution $(\n,u ), $
defined up to a positive time $T_0 $ which may depend on
$\inf\limits_{x\in \O}\n_0(x), $  and satisfying  \eqref{2.2}   and \eqref{2.3}.
We set \be \la{ss1}T^* =\sup\left\{T\,\left|\, \sup_{0\le t\le T}\|(\n ,u )\|_{H^2}<\infty\right\}.\right.\ee Clearly, $T^*\ge T_0.$
  We claim that
 \be   T^*=\infty.\ee  Otherwise,    $T^*<\infty.$ Then,   we claim that  there
exists a positive constant $ \hat{C}$ which may  depend  on $T^* $
and $\inf\limits_{x\in \O}\n_0(x)$   such that, for all  $0< T<
 T^*,$  \be\la{y12}\ba \sup_{0\le t\le T}
\| \n \|_{H^2}   \le \hat{C},\ea \ee where and what follows, $\hat C$ denotes some generic positive constant  depending on $T^*$ but   independent  of $T.$ This together with \eqref{pa1}    contradicts \eqref{ss1}. The estimates \eqref{pa1},   \eqref{ab19},   and  \eqref{ab12} directly follow  from \eqref{2.2}, Lemma \ref{310}, and  Propositions    \ref{le5} and \ref{caq1}.

It remains to prove \eqref{y12}.  First, standard calculations together with \eqref{pa1} yield that   for any $T\in (0,T^*),$
\be\la{b5.4} \inf\limits_{(x,t)\in\O\times (0,T)} \n(x,t)\ge \inf\limits_{x\in\O}\n_0(x)\exp\left\{-\int_0^T\|\div u\|_{L^\infty}dt\right\}\ge \hat C^{-1}.\ee

We define \be\la{5.5} \sqrt{\n}\dot u(x,t=0)=\n_0^{-1/2}\left( \mu \lap u_0+\na((\mu+\lambda(\n_0))\div u_0)-\na P(\n_0)\right).\ee
Integrating \eqref{b17} with respect to $t$ over $(0,T)$ together with \eqref{2.1}, \eqref{b3.56}, and \eqref{5.5}  yields \be \la{5.6}\sup\limits_{0\le t\le T}\int \n |\dot u|^2dx+\int_0^T\|\na\dot u\|_{L^2}^2dt\le \hat C.\ee
This combining with   \eqref{a39},   \eqref{4.19},  and \eqref{pa1} leads to
\be\la{a5.7}\ba & \sup\limits_{0\le t\le T}\left(\|\na^2u\|_{L^2}+\|\na G\|_{L^2}+\|\na \o\|_{L^2}\right)+\int_0^T(\|\na^2G\|_{L^2}^2+\|\na^2\o\|_{L^2}^2)dt\\ &\le \hat C \sup\limits_{0\le t\le T}\|\n\dot u\|_{L^2}+ \hat C \int_0^T\|\na(\n\dot u)\|_{L^2}^2dt\\ &\le \hat C + \hat C  \int_0^T\left(\|\na\n\|_{L^q}^2\|\dot u \|_{L^{2q/(q-2)}}^2+\|\na\dot u\|_{L^2}^2\right)dt\\ &\le \hat C + \hat C \int_0^T \| \dot u\|_{H^1}^2 dt\le \hat C,\ea\ee where in the last inequality, we have used \eqref{2.5} and \eqref{5.6}.

Next, operating $\na$ to \eqref{4.20} and  multiplying the resulting equality by $\na\Phi^i,$ we obtain after integration by parts and using \eqref{b5.4} and  \eqref{a5.7} that
\be\la{5.4}\ba \frac{\rm d}{{\rm d}t}\|\na\Phi\|_{L^2}   &\le \hat C(1+\|\na u\|_{L^\infty})\left(1+\|\na \Phi\|_{L^2}+\|\na\n\|_{L^4}^2+\|\na^2 \n\|_{L^2} \right)\\&\quad+\hat C\||\na\n||\na^2 u|\|_{L^2} + \hat C\||\na\n| |\na G|\|_{L^2}+ \hat C \|\na^2 G\|_{L^2}     .\ea\ee
Note that \eqref{lp} and  \eqref{pa1}   lead to\bnn\ba   \|\na^2\n\|_{L^2}+ \|\na \n\|_{L^4}^2 &\le \hat C\|\na \Phi\|_{L^2}+ \hat C\|\na \n\|_{L^4}^2\\ &\le \hat C\|\na \Phi\|_{L^2}+ \hat C\|\na \n\|_{L^{\min\{4,q\}}}^{\min\{4,q\}/2}\|\na^2 \n\|_{L^2}^{(4-\min\{4,q\})/2} \\ &\le  \hat  C\|\na \Phi\|_{L^2}+ \frac12\|\na^2 \n\|_{L^2} + \hat C,\ea\enn which together with  \eqref{b5.4} yields that \be \la{a509} \ba \|\na^2P(\n)\|_{L^2}+\|\na^2\n\|_{L^2}+\|\na \n\|_{L^4}^2 & \le  \hat  C\|\na^2\n\|_{L^2}+ \hat  C\|\na \n\|_{L^4}^2\\& \le  \hat  C\|\na \Phi\|_{L^2} + \hat C. \ea\ee

Then, on one hand, it follows from
Holder's inequality, \eqref{lp}, and \eqref{pa1} that
\be\la{a510}\ba &\||\na\n||\na^2 u|\|_{L^2} +  \||\na\n| |\na G|\|_{L^2}\\ &\le \hat C\|\na \n\|_{L^q}\|\na^2u\|_{L^2}^{1-2/q} \|\na^3u\|_{L^2}^{2/q}+\hat C\|\na \n\|_{L^q}\|\na G\|_{L^2}^{1-2/q} \|\na^2G\|_{L^2}^{2/q}\\ &\le \hat C(\ve)+\ve \|\na^3u\|_{L^2} +  \hat C  \|\na^2 G\|_{L^2}.\ea\ee
On the other hand, the $L^2$-estimate of elliptic system leads to
\be\la{a511}\ba \|\na^3u\|_{L^2}&=\|\na \lap u\|_{L^2}\\&\le  \hat C  \|\na^2\div u\|_{L^2}+
 \hat C  \|\na^2\o\|_{L^2}\\&\le  \hat C  \|\na^2((2\mu+\lambda(\n))\div u)\|_{L^2}+\hat C  \||\na\n|| \na^2 u |\|_{L^2}\\&\quad+\hat C  \||\na^2  \n ||\na u| \|_{L^2}+\hat C  \||\na   \n |^2|\na u| \|_{L^2}+
 \hat C  \|\na^2\o\|_{L^2} \\&\le  \hat C  \|\na^2G\|_{L^2}+\hat C\|\na^2\n\|_{L^2}+\hat C\|\na \n\|_{L^4}^2+
 \hat C  \|\na^2\o\|_{L^2}\\&\quad+\hat C  \||\na\n|| \na^2 u |\|_{L^2} +\hat C \left( \| \na^2  \n   \|_{L^2}+ \| \na \n   \|^2_{L^4}\right)\|\na u\|_{L^\infty}. \ea\ee
Substituting \eqref{a511} and  \eqref{a509} into  \eqref{a510} leads to
\bnn\la{a512}\ba & \||\na\n||\na^2 u|\|_{L^2} +  \||\na\n| |\na G|\|_{L^2}\\ &\le  \hat C  \|\na^2G\|_{L^2} +
 \hat C  \|\na^2\o\|_{L^2}+\hat C(1+\|\na
\Phi\|_{L^2} )(1+\|\na u\|_{L^\infty}), \ea\enn
which together with  \eqref{5.4}  and \eqref{a509} gives
\bnn\la{a5.4}\ba \frac{\rm d}{{\rm d}t}\|\na\Phi\|_{L^2}   &\le \hat C(1+\|\na u\|_{L^\infty})\left(\|\na \Phi\|_{L^2}+1\right)+ \hat C  \|\na^2G\|_{L^2} +
 \hat C  \|\na^2\o\|_{L^2}  .\ea\enn
This combining with  \eqref{a5.7} and Gronwall's inequality yields
 \bnn \sup\limits_{0\le t\le T}\|\na\Phi\|_{L^2}\le \hat C,\enn
  which together with \eqref{a509} implies \eqref{y12}. The proof of Proposition \ref{pro2} is finished.

{\it Proof of Theorem   \ref{t2}.}
 Let $(\n_0,m_0)$ satisfying  \eqref{1.9}  be the initial data as described in
Theorem \ref{t2}.  For constant   $ \de\in (0,1),$ we define
\be \la{5.15}\n_0^\de\triangleq   j_\de*\n_0  +\de\ge \de>0 ,\quad
 u_0^\de\triangleq   j_\de*u_0 ,\quad m_0^\de=\n_0^\de u_0^\de,\ee where  $j_\de$
is the standard mollifying kernel of width $\de.$ Hence, we have $\n_0^\de,u_0^\de \in H^\infty,$ and \bnn \la{5.16}\lim\limits_{\de\rightarrow 0}\left(\|\n_0^\de-\n_0\|_{ W^{1,q}}+\|u_0^\de-u_0\|_{H^1}\right)=0.\enn
Proposition \ref{pro2}  thus yields that the problem \eqref{n1}-\eqref{n4} with $(\n_0,m_0)$ being replaced by $(\n_0^\de,m_0^\de )$  has a unique global strong solution  $(\n^\de,u^\de)$  satisfying \eqref{pa1} for any $T>0$ and for  some $C$ independent of $\de.$ Moreover,  if \eqref{1.11} holds, there exists some positive constant $C$ depending only on $  \mu,  \beta, \gamma, \|\n_0\|_{L^\infty},$ and $\|u_0\|_{H^1} $ such that  $(\n^\de,u^\de)$  satisfies \eqref{ab12} and \eqref{ab19}. Letting $\de\rightarrow 0,$ standard   arguments (see \cite{Mik,Ka,hli1,lzz}) thus show that  the problem  \eqref{n1}-\eqref{n4}  has a global strong solution $(\n,u)$ satisfying the properties listed in Theorem \ref{t2} except \eqref{1.13} and the uniqueness of $(\n,u)$  satisfying \eqref{1.10}. Moreover,  $(\n ,u )$  satisfies \eqref{ab12} and \eqref{ab19} for some positive constant $C$ depending only on $  \mu,  \beta, \gamma, \|\n_0\|_{L^\infty},$ and $\|u_0\|_{H^1} $ provided \eqref{1.11} holds.

Since the uniqueness of $(\n,u)$  satisfying \eqref{1.10} is similar to that of Germain \cite{Ge} and \eqref{1.13} will be proved  in Theorem  \ref{t1}, we   finish the proof of Theorem  \ref{t2}.

{\it Proof  of  Theorem  \ref{t1}.}
 Let $(\n_0,m_0)$ satisfying  \eqref{v1.9}  be the initial data as described in
Theorem \ref{t1}.  For constant   $ \de\in (0,1),$ let $(\n_0^\de,u_0^\de)$  be as in \eqref{5.15}.  Hence, we have $\n_0^\de,u_0^\de \in H^\infty,$ and for any $p>1,$ \bnn \lim\limits_{\de\rightarrow 0}\left(\|\n_0^\de-\n_0\|_{ L^p}+\|u_0^\de-u_0\|_{H^1}\right)=0.\enn Moreover, \bnn\n_0^\de  \rightharpoonup \n_0 \mbox{ weakly * in }  L^\infty,\, \mbox{  as }\de\rightarrow 0.\enn
Proposition \ref{pro2}  thus yields that the problem \eqref{n1}-\eqref{n4} with $(\n_0,m_0)$ being replaced by $(\n_0^\de,\n_0^\de u_0^\de )$  has a unique global strong solution  $(\n^\de,u^\de)$  satisfying \eqref{b3.56}, \eqref{b19},   and  \eqref{4a2}, for any $T>0$ and for  some $C$ independent of $\de.$ Moreover,  if \eqref{1.11} holds, there exists some positive constant $C$ depending only on $  \mu,  \beta, \gamma, \|\n_0\|_{L^\infty},$ and $\|u_0\|_{H^1} $ such that  $(\n^\de,u^\de)$  satisfies \eqref{ab12},  \eqref{var31}, and \eqref{ab19}.

We modify the compactness  arguments  in \cite{Ka,Mik} to obtain the compactness results of $(\nd,\ud).$

First, it follows from \eqref{b3.56} and \eqref{4a2} that \bnn \la{5.30}\sup\limits_{0\le t\le T}\|\ud\|_{H^1}+\int_0^Tt\|\ud_t\|_{L^2}^2dt \le C,\enn
which together with the Aubin-Lions lemma gives that, up to a subsequence,
\bnn \la{5.31}\begin{cases}\ud\rightharpoonup u \mbox{ weakly * in }  L^\infty(0,T;H^1), \\ \ud\rightarrow u  \mbox{ strongly in } C([\tau,T];L^p), \end{cases}\enn for any $\tau\in (0,T)$ and $p\in [1,\infty).$

Next, let $A^\de\triangleq(2\mu+\lambda(\nd))\div \ud-P(\nd).$ One thus deduces from \eqref{ra3}, \eqref{rb3},  \eqref{b3.56}, and  \eqref{4a2} that
\be\la{5.32} \int_0^T\left( \|A^\de\|_{L^\infty}^{4/3}+ \|\o^\de\|_{H^1}^2+\|A^\de\|_{H^1}^2+t\|\o^\de_t\|_{L^2}^2+t\|A^\de_t\|_{L^2}^2
\right)dt\le C,\ee which implies that, up to a subsequence,
\be \la{5.33}\begin{cases}A^\de\rightharpoonup  A \mbox{ weakly * in }  L^{4/3}(0,T;L^\infty), \\ \o^\de \rightarrow \o=\curl u ,\,\, A^\de \rightarrow  A\mbox{ strongly  in } L^2 (\tau,T ;L^p),\end{cases} \ee for any $\tau\in (0,T)$ and $p\in [1,\infty).$

 Next, to obtain  the strong limit of $\nd,$ we deduce from \eqref{b3.56} that,   up to a subsequence, \bnn  \la{5.34} \nd\rightharpoonup \n \mbox{ weakly * in }  L^\infty(0,T;L^\infty).\enn
Let $f(s)$ be an arbitrary continuous function on $[0,C]$ with $C$ as in \eqref{b3.56}. Then, we have that,   up to a
subsequence, $f(\nd)$ converges weakly $*$  in $L^\infty(0,T;L^\infty).$ Denote the weak-$*$  limit by $\ol {f(\n)}:$
 \bnn \la{5.35}f(\nd)\rightharpoonup \ol {f(\n)} \mbox{ weakly * in }L^\infty(0,T;L^\infty). \enn
Noticing that,
\bnn \la{5.36}\div \ud =\phi(\nd)A^\de+\phi(\nd)P(\nd),\enn
with $\phi(s)\triangleq {1}/({2\mu+\lambda(s)}),$ we have  \bnn \la{5.37}\div u =\ol{\phi(\n)} A +\ol{\phi(\n) P(\n)} ,\quad \mbox{ a.e. in }\O\times (0,T). \enn

From \eqref{n1}, we obtain
\bnn  \la{5.39}(\ol {\n^2})_t+\div (\ol {\n^2}u)+A\ol{\n^2\phi(\n)}+\ol{\n^2\phi (\n) P(\n)}=0,\mbox{ in } \mathcal{D}'(\O\times (0,\infty)).
\enn
Using \cite[Lemma 2.3]{L2}, we get by standard arguments that
\bnn \la{5.40}( \n^2)_t+\div ( \n^2 u)+A\n^2\ol{\phi(\n)}+\n^2\ol{\phi(\n) P(\n)}=0, \mbox{ in } \mathcal{D}'(\O\times (0,\infty)).\enn
Thus, for $\Phi\triangleq\ol{\n^2}-  {\n^2}\ge 0, $ we have
\be  \la{5.41}\begin{cases}\Phi_t+\div(\Phi u)+A(\ol{\n^2 \phi(\n)}-\n^2\phi(\n))+A\n^2(\phi(\n)-\ol{\phi(\n)})+\ol{\n^2 \phi(\n) P(\n)}\\\quad- {\n^2 \phi(\n) P(\n)}+\n^2\left({ \phi(\n) P(\n)}-  \ol{ \phi(\n) P(\n)}\right)=0,\,\, \mbox{ in } \mathcal{D}'(\O\times (0,\infty)),\\ \Phi(x,t=0)=0,\quad \mbox{ a.e. }x\in \O.   \end{cases}\ee
By writing $(\nd)^2-\n^2=2\n(\nd-\n)+(\nd-\n)^2,$ we see that, up to a subsequence, \bnn  \la{5.43} \ol \lim_{\de\rightarrow 0}\|\nd-\n\|_{L^2}^2(t)\le  \int\Phi(x,t) dx,\quad \mbox{ for }t>0.\enn
   Also, for any $f(s)\in C^2([0,C])$ and $h(x)\in L^\infty(\O),$  noticing that \bnn \la{5.42} \ba f(\nd)-f(\n)&=f'(\n)(\nd-\n)+\int_0^1\te\int_0^1 f''(\n+\te\al (\nd-\n))d\al d\te(\nd-\n)^2 ,\ea\enn  we deduce from \eqref{b3.56} that \bnn \la{5.44} \left|\int h(x)(\ol{f(\n)}-f(\n))dx\right|\le M\|h\|_{L^\infty}\int \Phi dx,\enn for some constant $M>0.$
 In particular, noticing that \bnn f_1(s)\triangleq s^2 \phi(s)\in C^2([0,C]), \quad f_2(s)\triangleq s^2\phi(s)P(s)\in  C^2([0,C]), \enn we have \be  \la{5.45}\left|\int A\left(\ol{\n^2\phi(\n)}-\n^2\phi(\n)\right)dx\right|\le M\|A\|_{L^\infty}\int\Phi dx,\ee
 and \be  \la{5.46}\left|\int  (\ol{\n^2\phi(\n)P(\n)}-\n^2\phi(\n)P(\n))dx\right|\le M \int\Phi dx.\ee

 Let  $g(s)\in C^1([0,C])\cap C^2((0,C])$  be such that  for any $\n,s\in [0,C],$  \be \la{5.49}\left|\n^2 \int_0^1\te\int_0^1 g''(\n+\te\al (s-\n))d\al d\te\right|\le M.\ee
Note that
\bnn\ba &\n^2 (g(\nd)-g(\n))-\n^2  g'(\n)(\nd-\n )\\ &= \n^2 \int_0^1\te\int_0^1 g''(\n+\te\al (\nd-\n))d\al d\te (\nd-\n)^2,\ea\enn which together with \eqref{5.49} yields that for any   $h(x)\in L^\infty(\O),$ \be \la{5.50} \left|\int h(x)\n^2(\ol{ g(\n)}-  g(\n))dx\right|\le M\|h\|_{L^\infty}\int \Phi dx.\ee

Let $g_1(s)\triangleq\phi(s)$ and $g_2(s)\triangleq\phi(s)P(s).$ Since $g_i\in C^1([0,C])\cap C^2((0,C]) (i=1,2)$ satisfy \eqref{5.49},   from \eqref{5.50} we obtain that\be \la{5.51} \left|\int A\n^2(\ol{ \phi(\n)}-  \phi(\n))dx\right|\le M\|A\|_{L^\infty}\int \Phi dx,\ee and that \be\la{5.52}  \left|\int \n^2(\ol{ \phi(\n)P(\n)}-  \phi(\n)P(\n))dx\right|\le M \int \Phi dx.\ee

Substituting \eqref{5.45}, \eqref{5.46}, \eqref{5.51},  and \eqref{5.52} into \eqref{5.41}, after using Gronwall's inequality and \eqref{5.32}, we arrive at
\bnn \Phi=0 \mbox{ a.e.  in }\O\times (0,T),\enn which gives that, up to a subsequence, \be \la{5.54}\nd\rightarrow \n \mbox{ strongly in }L^p(\O\times (0,T)),\ee
for any $p\in [1,\infty).$ This combining with \eqref{5.33} implies that,  up to a subsequence, \be \la{l56} G^\de\rightarrow G= (2\mu+\lambda(\n))\div u+P(\n)-\bar P,\,\, \mbox{ strongly  in } L^2(\O\times (\tau,T)),\ee for any $\tau\in (0,T),$ where and what follows, $\ol f$ denotes the mean value of $f$ over $\O$ as in \eqref{1.6}.

Standard arguments thus show that the limit $(\n,u)$ is a global weak solution of \eqref{n1}-\eqref{n4}.

To finish the proof of Theorem \ref{t1}, it only remains to prove   \eqref{1.13}.

 First,
  it follows from \eqref{ab12} that
\bnn\la{za4} \int_1^\infty  \|P(\nd)-\ol {P(\nd)}\|_{L^2}^2 dt   \le C\int_1^\infty\left( \|G^\de\|_{L^2}^2+ (A_3^\de(t))^2\right)dt \le C,\enn
 which combining with \eqref{5.54}, \eqref{5.33},   \eqref{l56}, and  \eqref{ab12}  gives \be\la{zz5} \int_1^\infty  \left(\|P(\n )-\ol {P }\|_{L^2}^2+\|G\|_{L^2}^2+\|\o\|_{L^2}^2\right) dt    \le C.\ee
Simple calculations lead to
  \bnn  \ba &\frac{d}{dt}( \|P(\nd)-\ol {P(\nd)}\|_{L^2}^2 )
   \\&=2  \int(P(\nd)-\ol {P(\nd)})(P(\nd)-\ol {P(\nd)})_tdx \\&=-2  \int(P(\nd)-\ol {P(\nd)})(\ud\cdot \na(P(\nd)-\ol {P(\nd)})+\nd P'(\nd)\div \ud )dx\\&\quad +2  \int(P(\nd)-\ol {P(\nd)})dx \int (\nd P'(\nd)-P(\nd))\div \ud dx \\&\le C \|P(\nd)-\ol {P(\nd)}\|_{L^2}^2 +C\|\na \ud\|_{L^2}^2\\&\le C \|P(\nd)-\ol {P(\nd)}\|_{L^2}^2 +C\|G^\de\|_{L^2}^2+C\|\o^\de\|_{L^2}^2, \ea \enn
  which gives that, for any $s,t\in [N,N+1],$ \be\la{5.66} \ba &\|P(\nd)-\ol {P(\nd)}\|_{L^2}^2(t)-\|P(\nd)-\ol {P(\nd)}\|_{L^2}^2(s) \\ &\le   C\int_N^{N+1}\left( \|P(\nd)-\ol {P(\nd)}\|_{L^2}^2+\|G^\de\|_{L^2}^2+ \|\o^\de\|_{L^2}^2   \right)dt.\ea\ee  Integrating \eqref{5.66} with respect to $s$ over $(N,N+1)$ yields that \bnn \ba &\sup\limits_{N\le t\le  N+1}\|P(\nd)-\ol {P(\nd)}\|_{L^2}^2(t)\\ &\le C\int_N^{N+1}\left( \|P(\nd)-\ol {P(\nd)}\|_{L^2}^2+\|G^\de\|_{L^2}^2+ \|\o^\de\|_{L^2}^2 \right)dt.\ea\enn
    From \eqref{5.54},  \eqref{l56}, and \eqref{5.33}, we have
\bnn \ba  \sup\limits_{N\le t\le  N+1}\|P(\n )-\ol {P}\|_{L^2}^2(t) \le C\int_N^{N+1}\left( \|P(\n )-\ol {P }\|_{L^2}^2+\|G \|_{L^2}^2+ \|\o \|_{L^2}^2 \right)dt.\ea\enn
  Letting $N\rightarrow \infty,$ this combining with \eqref{zz5} yields that
  \be \la{5.69}\lim\limits_{t\rightarrow 0}\|P(\n )-\ol {P }\|_{L^2}^2(t)=0.\ee

Next, standard arguments together with \cite[Lemma 2.3]{L2} and  $(\ref{n1})_1$ yield  that $P(\n)$ satisfies
\bnn (P(\n))_t+\div(P(\n)u)+(\ga-1)P(\n)\div u=0,\quad \mbox{ in }\mathcal{D}'(\O\times(0,\infty)),\enn
which gives that
  \bnn \la{za8}\ba \int_1^\infty \left|\frac{d}{dt}\ol{ P  }\right|dt&\le C\int_1^\infty\left|\int (P-\bar P)\div  u dx\right|dt\\ &\le C\int_1^\infty\left( \|P-\bar P\|_{L^2}^2+\|\na u\|_{L^2}^2 \right)dt \le C,\ea\enn  due to \eqref{ab12}.
  Hence, there exists some positive constant $\n_s$ such that
\bnn \la{za9}\lim\limits_{t\rightarrow \infty}  \bar P (t)=\n_s^\ga,\enn
due to $0<\bar\n_0^\ga\le \bar P\le C.$ This combining with \eqref{5.69} and \eqref{ab12} yields that\bnn  \la{za10}\lim\limits_{t\rightarrow \infty}\|\n-\n_s\|_{L^p}(t)=0,\enn for any $p\in [1,\infty).$  Thus, \eqref{c3.10} gives \be \la{za11} \lim\limits_{t\rightarrow \infty}\|\n-\bar\n_0\|_{L^p}(t)=0,\ee for any $p\in [1,\infty).$

Finally, similar to \eqref{5.66}, from \eqref{var31} and \eqref{ab12}, we have \bnn \ba (A_1^\de(t)  )^2
\le (A_1^\de (s))^2 +C\int_N^{N+1} (A_3^\de(t))^2  dt,\ea\enn for any $s,t\in [N,N+1].$ This gives \bnn \ba  \sup\limits_{N\le t\le N+1}(A_1^\de(t)  )^2 & \le   C\int_N^{N+1} \left((A_1^\de (t))^2+(A_3^\de(t))^2 \right) dt\\&\le C\int_N^{N+1}\left( \|P(\nd )-\ol {P(\nd )}\|_{L^2}^2+\|G^\de\|_{L^2}^2+ \|\o^\de\|_{L^2}^2 \right)dt,\ea\enn
which together with \eqref{5.33}, \eqref{l56},  \eqref{zz5}, and the fact that
$$ \|G^\de\|_{L^2}^2+\|\o^\de\|_{L^2}^2  \le C(A_1^\de(t)  )^2,$$
 leads to
  \bnn\lim\limits_{t\rightarrow 0}(\|G \|_{L^2}^2+\|\o \|_{L^2}^2 ) (t)=0.\enn
Because of \eqref{5.69}, this shows
\bnn\lim\limits_{t\rightarrow \infty}\|\na u\|_{L^2}\le C
  \lim\limits_{t\rightarrow \infty}(\|G \|_{L^2} +\|\o \|_{L^2}+ \|P-\bar P\|_{L^2})(t) =0,\enn
which combining with \eqref{ab19} and \eqref{za11}  directly yields
  \eqref{1.13}. The proof  of  Theorem  \ref{t1} is completed.

{\it Proof  of  Theorem  \ref{th2}.}  Since the proof of  Theorem  \ref{th2} is similar to that of  \cite[Theorem 1.2]{lx}, we omit it here.

\end{document}